\documentclass[11pt,openbib]{article}
\usepackage{amsmath}
\usepackage{amssymb,latexsym,theorem}

\theoremstyle{plain}
\theorembodyfont{\normalfont}
\newtheorem{theorem}{Theorem}
\makeatletter
\renewcommand\tiny{\@setfontsize\tiny\@viiipt\@xpt}
\makeatother
\theorembodyfont{\rmfamily}
\newtheorem{thm}[theorem]{Theorem}

\newtheorem{prop}[theorem]{Proposition}

\newtheorem{lem}[theorem]{Lemma}

\theoremstyle{definition}
\newtheorem{example}[theorem]{Example}

\theoremstyle{remark}
\newtheorem{remark}[theorem]{Remark}
\newtheorem{fact}[theorem]{Fact}

\newcommand{\R}{\mathbb{R}}

\newcommand{\Aff}{\mathrm{Aff}\!}
\newcommand{\aff}{\mathrm{aff}\!}
\newcommand{\gl}{\mathop{\mathrm{gl}}\nolimits}
\newcommand{\Gl}{\mathop{\mathrm{Gl}}\nolimits}

\newcommand{\oo}{\mathop{\mathrm{\! \, o}}\nolimits}
\newcommand{\Oo}{\mathop{\mathrm{\! \, O}}\nolimits}

\newcommand{\Ad }{\mathop{\mathrm{Ad}}\nolimits}

\newcommand{\comp}{\raisebox{0pt}{$\scriptstyle\circ \, $}}
\newcommand{\setrule}{\, \rule[-4pt]{.5pt}{13pt}\, }

\newcommand{\verysmallrowspace}{\rule{0pt}{10pt}}
\newcommand{\smallrowspace}{\rule{0pt}{14pt}}
\newcommand{\rowspace}{\rule{0pt}{16pt}}

\newcommand{\ttfrac}[2]{\mbox{${\scriptstyle \frac{{#1}}{{#2}}}$}}
\newcommand{\onehalf}{\mbox{$\frac{\scriptstyle 1}{\scriptstyle 2}\,$}}
\newcommand{\spann}{\mathop{\rm span}\nolimits}

\newcommand{\dbydt}{\mbox{${\scriptstyle \frac{d}{dt}}
\rule[-6pt]{.5pt}{15pt} \raisebox{-6pt}{$\, {\scriptscriptstyle t=0}$}$}}

\begin{document}

\begin{center}
{\Large \bf Adjoint and coadjoint orbits
\rule{10pt}{0pt} \\ of the Poincar\'{e} group} \\
\mbox{} \vspace{.05in} \\
Richard Cushman and Wilberd van der Kallen\footnotemark
\end{center}
\footnotetext{Mathematics Institute, University of Utrecht, 3508TA Utrecht,
The Netherlands }
\addtocounter{footnote}{1}
\footnotetext{This is a revised version of \cite{cushman-vdk}.}
\addtocounter{footnote}{1}
\footnotetext{printed: \today}

\begin{abstract}
In this paper we give an effective method for finding
a unique representative of each orbit of the adjoint and coadjoint action
of the real affine orthogonal group on its Lie algebra. In both cases
there are orbits which have a modulus that is different from the
usual invariants for orthogonal groups. We find an 
unexplained bijection between adjoint and coadjoint orbits. As a
special case, we classify the adjoint and coadjoint orbits of the
Poincar\'{e} group.
\end{abstract}

\section{Introduction}\label{intro}

Let $(\widetilde{V}, \widetilde{\gamma})$ be an $n$-dimensional
real vector space with a nondegenerate inner product $\widetilde{\gamma}$.
The set $\Oo (\widetilde{V}, \widetilde{\gamma })$ of real linear maps
$B$ of $\widetilde{V}$ into itself, which preserve $\widetilde{\gamma }$,
that is, $\widetilde{\gamma }(Bv,Bw) = \widetilde{\gamma }(v,w)$ for
every $v,w \in \widetilde{V}$, is a Lie group called the
\emph{orthogonal group}. Its Lie algebra $\oo (\widetilde{V},
\widetilde{\gamma})$ consists of real linear maps $\xi $ of $\widetilde{V}$
into itself such that $\widetilde{\gamma }(\xi v, w) +
\widetilde{\gamma }(v, \xi w) =0$ for every $v,w \in \widetilde{V}$.
For $\xi , \, \eta \in \widetilde{V}$
the Lie bracket on $\oo (\widetilde{V}, \widetilde{\gamma})$ is
$[ \xi , \eta ] = \xi \comp \eta - \eta \comp \xi $, where
$\comp $ is the composition of linear maps. The \emph{affine orthogonal
group} $\Aff \Oo (\widetilde{V}, \widetilde{\gamma }) =
\Oo (\widetilde{V}, \widetilde{\gamma }) \ltimes
\widetilde{V}$ is the set of real affine orthogonal maps of $(\widetilde{V},
\widetilde{\gamma })$ into itself. More precisely, it is the set
$\Oo (\widetilde{V}, \widetilde{\gamma }) \times \widetilde{V}$
with group multiplication $(B_1, v_1) \cdot (B_2, v_2) = 
(B_1B_2, B_1v_2 +v_1)$, which is the composition of affine linear maps. 
The affine orthogonal group is a Lie group. Its Lie algebra
$\aff \oo (\widetilde{V}, \widetilde{\gamma }) =
\oo (\widetilde{V}, \widetilde{\gamma}) \times
\widetilde{V}$ has
Lie bracket $\left[ ({\xi }_1, v_1),({\xi }_2, v_2) \right] = 
( [ {\xi }_1, {\xi }_2],
{\xi }_1v_2 - {\xi }_2v_1 ) $, where ${\xi }_1, {\xi }_2$ lie in 
$\oo (\widetilde{V}, \widetilde{\gamma })$. The \emph{adjoint action} of
the affine orthogonal group on its Lie algebra is defined by
\begin{displaymath}
\begin{array}{l}
\varphi :(\Oo (\widetilde{V}, \widetilde{\gamma }) \ltimes \widetilde{V})
\times (\oo (\widetilde{V}, \widetilde{\gamma }) \times \widetilde{V})
\rightarrow \oo (\widetilde{V}, \widetilde{\gamma }) \times 
\widetilde{V}: \\
\rowspace \hspace{.25in} ((B, v), (\xi , w)) \mapsto
(B,v) \cdot (\xi , w) \cdot {(B,v)}^{-1},
\end{array}
\end{displaymath}
where $\cdot $ is composition of affine linear maps. A straightforward
calculation shows that ${\varphi }((B,v), (\xi ,w)) = 
(B\xi B^{-1}, -B\xi B^{-1}v + Bw)$.

One of the goals of this paper is to classify the orbits of the adjoint
action of the affine orthogonal group. In particular, we find a unique
representative (= normal form) for each orbit. The basic technique leans
heavily on the idea of an indecomposable type introduced by
Burgoyne and Cushman \cite{burgoyne-cushman} to find normal forms
for the adjoint action  of any real form of a nonexceptional 
Lie group.\footnote{We recall the idea of an indecomposable 
type for the special case of the Lie algebra
$\gl (\widetilde{V})$ of the Lie group $\Gl (\widetilde{V})$ of
invertible real linear maps of $\widetilde{V}$ into itself.
Let $\widetilde{\xi } \in \gl (\widetilde{V})$. Consider the
\emph{pair} $(\widetilde{\xi }, \widetilde{V})$. On the collection
of all pairs we say that two pairs $(\widetilde{\xi }, \widetilde{V})$
and $({\widetilde{\xi }}', {\widetilde{V}}')$ are \emph{equivalent}
if there is a bijective real linear map
$P:\widetilde{V} \rightarrow \widetilde{V}'$
such that $P \widetilde{\xi } = {\widetilde{\xi }}' P$. Note that
$P$ defines an isomorphism of $\Gl (\widetilde{V})$ with
$\Gl ({\widetilde{V}}')$. Clearly being equivalent is an equivalence
relation on the collection of pairs. We call an equivalence class a
\emph{type}. Let $\Delta $ be the
type represented by the pair $(\widetilde{\xi }, \widetilde{V})$.
Suppose that $\widetilde{V} = {\widetilde{W}}_1 \oplus
{\widetilde{W}}_2$, where ${\widetilde{W}}_i$ are proper,
$\widetilde{\xi }$-invariant subspaces, then
$\widetilde{\xi }|{\widetilde{W}}_i \in \gl ({\widetilde{W}}_i)$.
Let ${\Delta }_i$ be the type represented by
$(\xi |{\widetilde{W}}_i, {\widetilde{W}}_i)$. Then $\Delta $ is
the \emph{sum} of ${\Delta }_1$ and ${\Delta }_2$, which we write as
$\Delta = {\Delta }_1 + {\Delta }_2$. This sum is well defined.
We say the the type $\Delta $ is \emph{indecomposable} if it can
not be written as the sum of two or more types. An indecomposable
type is a Jordan block. The main theorem for classifying
the orbits of the adjoint action of $\Gl (\widetilde{V})$ on
$\gl (\widetilde{V})$ is: for every $\widetilde{\xi } \in
\gl (\widetilde{V})$, the type $\Delta $ represented by
$(\widetilde{\xi }, \widetilde{V})$ may be written as a sum of
indecomposable types, which is unique up to reordering of the
summands.  This is nothing but another formulation of the real
canonical form for real linear maps.} In this method the emphasis 
is not on subgroups and subvarieties, but rather on
vector spaces with quadratic forms. (Indeed we learn little about an orbit
as a variety. There is ample room for further work.)

Our aims are rather limited, but still we get results that seem to be new,
despite a widespread belief that all is known on this topic.
As explained in section \ref{section 1}
below, our affine orthogonal group may be viewed as a subgroup
of a slightly larger orthogonal group $O(V,K)$.
We find that the usual eigenvalue and Jordan invariants that 
classify the adjoint orbits of this ambient group $O(V,K)$ do 
not suffice to distinguish the orbits of the affine orthogonal 
group. That is why we have to invent a 
modulus, which parametrizes families of adjoint orbits, each family
being contained in a single orbit of $O(V,K)$.
In our classification of adjoint orbits we use the fact 
that we are working over the reals. 

Next let us turn to the classification of coadjoint orbits. Recall that Rawnsley \cite{rawnsley} has described how in principle 
one can classify the coadjoint orbits by reducing the problem to a 
similar problem for a subgroup known as the little subgroup. One 
should be careful though, because there is no canonical isomorphism between
the little subgroup as an actual subgroup and
your favorite incarnation of the isomorphism type of the little subgroup
as a Lie group. This matters 
because affine orthogonal groups are less rigid than ordinary orthogonal
groups. In particular, rescaling the vector part of an affine orthogonal 
group
gives an automorphism that is not inner. Thus performing the actual 
classification, as opposed to giving an in principle classification, 
needs some care. We do the classification in the style of 
Burgoyne and Cushman
\cite{burgoyne-cushman}, working with vector spaces instead 
of subgroups or subvarieties. Again we encounter an unfamiliar modulus.
Surprisingly, once we have found representatives of coadjoint orbits, we see
that there is a bijection between the  chosen
representatives for adjoint orbits and those employed for coadjoint orbits.
This bijection preserves ``dimension", ``index", ``modulus" and Jordan type.
We have no geometric explanation for it. \medskip 

We now give an overview of the contents of this paper. In
section \ref{section 1} we show that the affine orthogonal group 
is isomorphic to a larger orthogonal group, which leaves an isotropic
vector $v^{\comp }$ fixed. Throughout the remainder of the paper
we look only at this isotropy group. In section \ref{section 2} we adapt the
notion of an indecomposable type to the case at hand and show that
there is a distinguished indecomposable type containing
the vector $v^{\comp }$. In section \ref{section 3}
we classify these distinguished
indecomposable types and complete the classification of the
adjoint orbits of the affine orthogonal group. In section \ref{section 4} we
apply the above theory to find normal forms for the adjoint orbits
of the Poincar\'{e} group. In section \ref{section 5} we classify the 
coadjoint orbits of the affine orthogonal group and in section
\ref{section 6}
we specialize this to the coadjoint orbits of the Poincar\'{e} group. 

\section{Affine orthogonal group}\label{section 1}

In this section we show that the affine orthogonal group can be
realized as an isotropy subgroup of a larger orthogonal group. \medskip

Let $\widetilde{\gamma }$ be a nondegenerate inner product on
a real $n$-dimensional vector space $\widetilde{V}$. Let
$\{ e_1, \ldots \, , e_n \}$ be an orthonormal basis of
$\widetilde{V}$ such that the matrix of $\widetilde{\gamma }$ with
respect to this basis is $G = \mathrm{diag}(-I_m, I_p)$, where
$I_r$ is the $r \times r$ identity matrix. Let $\Oo (\widetilde{V}, G)$ be
the set of all linear maps $B$ of $\widetilde{V}$ into itself which
preserve $\widetilde{\gamma }$, that is, $\widetilde{\gamma }(Bv,Bw)
= \widetilde{\gamma }(v,w)$ for every $v,w \in \widetilde{V}$. Then
$\Oo (\widetilde{V}, G)$ is a Lie group which is isomorphic to
$\Oo (m,p)$. On $V = \R \times \widetilde{V} \times \R $ consider the
inner product $\gamma $ defined by $\gamma ((x,v,y), (x',v',y')) = 
\widetilde{\gamma }(v,v') + x'y + xy'$. With respect to the basis
$\mathfrak{e} = \{ e_0, e_1, \ldots \, , e_n, e_{n+1} \} $ of
$V$ the matrix of $\gamma $ is
\emph{standard}, that is, $K=${\tiny $\begin{pmatrix}
0 & 0 & 1 \\ 0 & G & 0 \\ 1 & 0 & 0 \end{pmatrix} $}.
Note that $e_{n+1}$ is a $K$-\emph{isotropic} vector of $(V,K)$, that is,
$K(e_{n+1},e_{n+1}) =0$. Let $\Oo (V, K)$ be the set of all
real linear maps $A$ of $V$ into itself
which preserve $\gamma $, that is, $\gamma (A(x,v,y), A(x',v',y')) =
\gamma ((x,v,y), (x',v',y'))$. \medskip

Now consider the isotropy subgroup
\begin{displaymath}
{\Oo (V,K)}_{e_{n+1}} = \{ A \in \Oo (V,K) \, \setrule \, Ae_{n+1} =
e_{n+1} \}
\end{displaymath}
of $\Oo (V,K)$. To give a more explicit description of
${\Oo (V,K)}_{e_{n+1}}$ let $A$ be an invertible real linear map 
of $V$ into itself such that $Ae_{n+1} = e_{n+1}$. Suppose that
the matrix of $A$ with respect to the basis $\mathfrak{e}$ is
{\tiny $\left( \begin{array}{ccc}
a & b^T & c \\ d & B & e \\ f & g^T & h \end{array} \right) $}. Then $A=$
{\tiny $\left( \begin{array}{ccc} a
& b^T & 0 \\
d & B & 0 \\ f& g^T & 1 \end{array} \right) $}, because $A$ leaves 
the vector
$e_{n+1}$ fixed. Now $A \in
\Oo (V,K)$ if and only if $K = A^TKA$, that is,
\begin{equation}
A = \mbox{{\tiny $\begin{pmatrix}
1 & 0 & 0 \\ d & B & 0 \\ -\onehalf d^TGd & -d^TGB & 1
\end{pmatrix} $},}
\label{eq-s1onestar}
\end{equation}
where $B^TGB = G$ and $d \in {\R }^n$. Thus $A \in {\Oo(V, K)}_{e_{n+1}}$
if and only if (\ref{eq-s1onestar}) holds. The group
${\Oo (V, K)}_{e_{n+1}}$ is isomorphic to the \emph{affine
orthogonal} group $\Aff \Oo (V, K)$, which is the semidirect product 
$\ltimes $
of $\Oo ({\R }^n, G)$ with ${\R }^n$, that is,
\begin{displaymath}
\Oo ({\R }^n, G) \ltimes  {\R }^n =
\left\{ \mbox{{ \tiny $\left( \begin{array}{cc}
1 & 0 \\ d & B \end{array} \right)$}} \in \Gl (n+1, \R ) \,
\setrule \, B^TGB =G, \, \, d \in {\R }^n \right\} .
\end{displaymath}
Explicitly, the isomorphism is given by
\begin{displaymath}
{\Oo (V, K)}_{e_{n+1}} \rightarrow \Oo ({\R }^n, G)
\ltimes  {\R }^n: \mbox{{  \tiny $\left( \begin{array}{ccc}
1 & 0 & 0 \\ d & B & 0 \\ -\onehalf d^TGd & -d^TGB & 1 \end{array}
\right) $}} \mapsto
\mbox{{  \footnotesize $\left( \begin{array}{cc} 1 & 0 \\
d & B \end{array} \right) $}}.
\end{displaymath}

We determine the Lie algebra ${\oo (V,K)}_{e_{n+1}}$ of
${\Oo (V, K)}_{e_{n+1}}$ as follows. Let  $v \in {\R }^n$ and
$X \in \oo (\widetilde{V}, G)$, that is, $X^TG + GX =0$. Then
\begin{displaymath}
t \mapsto \mbox{{\tiny $\left( \begin{array}{ccc}
1 & 0 & 0 \\ tv & \exp tX & 0 \\ -\onehalf (tv)^TG(tv) & -(tv)^TG 
\exp tX & 1
\end{array} \right)$ }} = Y_t
\end{displaymath}
is a curve in ${\Oo (V, K)}_{e_{n+1}}$ which passes through the
identity element at $t=0$.
Consequently, $\dbydt \hspace{-6pt} Y_t =${ \tiny $\begin{pmatrix}
0 & 0 & 0 \\ v & X & 0 \\ 0 & -v^TG & 0
\end{pmatrix} $} is an element of ${\oo (V,K)}_{e_{n+1}}$. 
The Lie bracket $[\, , \, ]$ on
${\oo (V,K)}_{e_{n+1}}$ is given by
\begin{displaymath}
\left[ \mbox{{ \tiny $\left( \begin{array}{ccc}
0 & 0 & 0 \\ v_1 & X_1 & 0 \\ 0 & -{v_1}^TG & 0
\end{array} \right) $}},
\mbox{{\tiny $\left( \begin{array}{ccc}
0 & 0 & 0 \\ v_2 & X_2 & 0 \\ 0 & -{v_2}^TG & 0
\end{array} \right) $}} \right]   {=} \mbox{{ \tiny $\left( 
\begin{array}{ccc}
0 & 0 & 0 \\ X_1v_2-X_2v_1 & [X_1,X_2] & 0 \\ 0 & -(X_1v_2-X_2v_1)^TG & 0
\end{array} \right) $},}
\end{displaymath}
where $[X_1,X_2]$ is the Lie bracket in $\oo (\widetilde{V}, G)$.

\section{Classification of adjoint orbits}\label{section 2}

To fix notation. Let $v^{\comp }$ be a nonzero \emph{isotropic} vector in
the real inner product space $(V, \gamma )$. Let
${\oo (V, \gamma )}_{v^{\comp }}$ be the Lie algebra
of the affine orthogonal group
\begin{displaymath}
{\Oo (V, \gamma )}_{v^{\comp }} = \{ A \in \Gl (V) \, \setrule
\, Av^{\comp } = v^{\comp } \, \, \mathrm{and} \, \, A^{\ast }\gamma 
= \gamma
\}.
\end{displaymath}
Then $Y \in {\oo (V, \gamma )}_{v^{\comp }}$ if and only if
\begin{displaymath}
Yv^{\comp } = 0 \, \, \mathrm{and} \, \, \gamma (Yv,w) +\gamma (v, Yw) =0,
\, \, \mbox{for all $v,w \in V$}.
\end{displaymath}

We begin our classification of the adjoint
orbits of the affine orthogonal group ${\Oo (V, \gamma )}_{v^{\comp }}$
on its Lie algebra ${\oo (V, \gamma )}_{v^{\comp }}$ by defining the notions
of indecomposable type and indecomposable distinguished type. First
we define the notion of a pair. Let $W$ be a $\gamma $-nondegenerate
real vector space. Our vector spaces are always finite dimensional.
If $Y \in \oo (W, \gamma )$
then $(Y, W;\gamma )$ is a \emph{pair}.\footnote{Our concept of
pair is the same as that of \cite{burgoyne-cushman}.} We say that the pairs
$(Y,W; \gamma )$ and $(Y', W' ;{\gamma }')$ are \emph{equivalent} if
there is a bijective real linear map $P:W \rightarrow W'$ such that
$PY = Y'P$ and $P^{\ast }{\gamma }' = \gamma $, that is,
${\gamma }'(Pv, Pw) = \gamma (v,w)$ for every $v,w \in W$. Clearly
being equivalent is an equivalence relation on the collection of pairs. An
equivalence class of pairs is a $\emph{type}$, which we denote by
$\Delta $. Given a type $\Delta $ with representative $(Y,W; \gamma)$
we define the \emph{dimension}, denoted $\dim \Delta $, of $\Delta $ by
$\dim W$ and the \emph{index}, denoted $\mathrm{ind}\, \Delta $, of
$\Delta $ by the number of negative eigenvalues of the Gram matrix
$(\gamma (v_i, v_j))$, where $\{ v_1, \dots \, , v_{\dim W} \} $ is a
basis of $W$. It is straightforward to check that neither of these
notions depends on the choice of representative of $\Delta $ or
on the choice of basis. Let $Y =S+N$ be the Jordan decomposition of
$Y$ into a semisimple linear map $S$ and a commuting nilpotent
linear map $N$, which lie in $\oo (W, \gamma )$. Because $S$ and $N$ are polynomials in $Y$ with
real coefficients and $Yv^{\comp }=0$, it follows that $Sv^{\comp } =
Nv^{\comp } = 0$. So $S,N \in {\oo (W, \gamma )}_{v^{\comp }}$. Let
$h$ be the unique nonnegative integer such that $N^hW \ne \{ 0 \}$ but
$N^{h+1}W =0$. We call $h$ the \emph{height} of the type $\Delta $ and
we denote it by $\mathrm{ht}(\Delta )$. It is evident that
$\mathrm{ht}(\Delta )$ does not depend of the choice of representative
of $\Delta $. We say that a type $\Delta $ with representative
$(Y, W;\gamma )$ is \emph{uniform} if $NW = \ker N^h|W$. Let
$(Y, W;\gamma )$ represent the type $\Delta $. Suppose that
$W = W_1 + W_2$, where $W_i$ are proper, $Y$-invariant subspaces,
which are $\gamma $-nondegenerate and $\gamma $ orthogonal. Then
$\Delta $  is the \emph{sum} of two types ${\Delta }_i$,
which are represented by $(Y|W_i, W_i; \gamma |W_i)$. We write 
$\Delta = {\Delta }_1+{\Delta }_2$. The type $\Delta $ is
\emph{indecomposable} if it can not be written as the sum of two
types. From \cite[prop.\ 3, p.343]{burgoyne-cushman} it follows that an
indecomposable type is uniform. So far the vector $v^{\comp }$ has
not played any role. Therefore the classification of indecomposable
types is given by results in \cite{burgoyne-cushman}. \medskip

We now define the notion of a triple, where the vector
$v^{\comp }$ plays an essential role. $(Y, W, v^{\comp };
\gamma )$ is a \emph{triple} if and only if the vector
$v^{\comp }$ is nonzero and $\gamma $-isotropic and for the linear map $Y$
in the pair $(Y,W; \gamma)$ we have $Yv^{\comp }=0$. We say that
the triple $(Y, W, v^{\comp };
\gamma )$ is a \emph{nilpotent triple} if $Y$ is
nilpotent. Two triples $(Y, W, v^{\comp };
\gamma )$ and $(Y', W', (v^{\comp })'; {\gamma }')$ are
\emph{equivalent} if there is a bijective real linear map
$P: W \rightarrow W'$ such that $Y'P = PY$, $P^{\ast }{\gamma }' =
\gamma $ and $Pv^{\comp } = (v^{\comp })'$. Clearly being equivalent
is an equivalence relation on the collection of triples. We call an 
equivalence
class of triples a \emph{distinguished type}, which we denote by
$\underline{\Delta }$. Let
$(Y, W, v^{\comp }; \gamma )$ represent the distinguished type
$\underline{\Delta }$. If $Y$ is nilpotent, then $\underline{\Delta }$ is
a \emph{nilpotent distinguished type}. Suppose that $W = W_1 \oplus W_2$,
where $W_i$ are proper, $Y$-invariant, $\gamma $-orthogonal,
$\gamma $-nondegenerate subspaces and $v^{\comp } \in W_1$. Then
$(Y|W_1, W_1,v^{\comp}; {\gamma }|W_1)$ is a triple whose
distinguished type we write ${\underline{\Delta }}_1$. Moreover, let
the pair $(Y|W_2, W_2; {\gamma }|W_2)$ represent the type
${\Delta }_2$. In this situation we say that the distinguished
type $\underline{\Delta }$ is the \emph{sum} of the distinguished
type  ${\underline{\Delta }}_1$ and the type ${\Delta }_2$ and
we write $\underline{\Delta } = {\underline{\Delta }}_1 +
{\Delta }_2$. If $\underline{\Delta }$ can not be written as the
sum of a distinguished type and a type, then we say that
$\underline{\Delta }$ is an \emph{indecomposable} distinguished type.
In other words, $(Y,W,v^{\comp };\gamma )$ represents an indecomposable
distinguished type if there is no proper, $\gamma $-nondegenerate,
$Y$-invariant subspace of $W$ which contains $v^{\comp }$. To
simplify notation from now on we usually drop the inner product
$\gamma $ in pairs and triples. \medskip

The first goal of this paper is to prove 

\begin{thm}\label{theorem1} 
Every distinguished type is
a sum of an indecomposable nilpotent distinguished type and a sum of
indecomposable types. This decomposition is unique up to a reordering 
of the summands. 
\end{thm} 

The proof of the theorem will require an understanding of
indecomposable nilpotent distinguished types.
Recall the indecomposable types have already been classified in
\cite{burgoyne-cushman}. The theorem solves the conjugacy 
class problem for the Lie algebra
${\oo (v, \gamma )}_{v^{\comp}}$.
Indeed distinguished types represented by triples
of the form $(Y,V,v^{\comp }; \gamma )$
correspond one to one with orbits of the adjoint action on
${\oo (v, \gamma )}_{v^{\comp}}$.\medskip

Before beginning the proof of theorem 1, we need some
additional concepts. Let $\underline{\Delta }$ be a distinguished type with
representative $(Y,W, v^{\comp})$. We say that $\underline{\Delta }$ has
\emph{distinguished height} $h$, if $h$ is the largest positive
integer for which there is a vector $w \in W$ such that
$Y^hw =v^{\comp }$. We denote the distinguished height of
$\underline{\Delta }$ by $\mathrm{dht}(\underline{\Delta })$. Because
the definition of distinguished height does not involve the
inner product $\gamma $ and $Yv^{\comp }=0$, there is a largest Jordan block
of the linear map $Y$ which contains the vector $v^{\comp }$. Moreover, 
it is of size $h+1$. Let
\begin{displaymath}
\mu (\underline{\Delta }) = \{ \gamma (w, v^{\comp }) \in \R \,
\setrule \, \mbox{for all $w \in W$ such that $Y^hw=v^{\comp }$} \} .
\end{displaymath}
We call $\mu (\underline{\Delta })$ the set of \emph{parameters} of the
distinguished type $\underline{\Delta }$. Below we will show that this set 
is a singleton. 

We prove  

\begin{lem}\label{Lemma 1.2} 
Suppose that
$\underline{\Delta } = {\underline{\Delta }}' + \Delta $. Then
$\mathrm{dht}(\underline{\Delta }) =
\mathrm{dht}({\underline{\Delta }}')$ and $\mu (\underline{\Delta }) =
\mu ({\underline{\Delta }}')$. 
\end{lem}

\noindent \textbf{Proof.} Suppose that $(Y,W, v^{\comp })$ is
a triple which represents the distinguished type $\underline{\Delta }$ and
that $W = W_1\oplus W_2$, where $W_i$ are proper, $Y$-invariant,
$\gamma $-orthogonal, $\gamma $-nondegenerate subspaces of $W$ with
$v^{\comp } \in W_1$. Say the triple $(Y|W_1, W_1, v^{\comp })$ represents a
distinguished type ${\underline{\Delta }}'$ and the pair
$(Y|W_2, W_2)$ represents the type $\Delta $. Suppose that
$\mathrm{dht}({\underline{\Delta }}') = h'$. Then there is a vector
$w' \in W_1$ such that $Y^{h'}w' = v^{\comp}$. Consequently,
$\mathrm{dht}(\underline{\Delta }) \ge h' $. Since
$\mathrm{dht}(\underline{\Delta })
= h$, there is a vector $w \in W$ such that $Y^hw= v^{\comp }$. But
$W = W_1\oplus W_2$. So we may write $w = w_1+w_2$ where $w_i \in W_i$.
Since $W_i$ are $Y$-invariant, we have $v^{\comp } =Y^hw_1 +Y^hw_2$ where
$Y^hw_i \in W_i$. By construction $v^{\comp } \in W_1$. Therefore
$Y^hw_1 = v^{\comp }$. Consequently $h \le \mathrm{dht}({\underline{\Delta }}')
= h'$. So $h = h'$. Note that $\dim \underline{\Delta } >
\dim {\underline{\Delta }}'$. \medskip

Since $W_1 \subseteq W$, it follows from the definition of the
set of parameters that $\mu ({\underline{\Delta }}')
\subseteq \mu (\underline{\Delta })$. Suppose
that there is a vector $w \in W$ with $Y^hw=v^{\comp }$ such that
$\gamma (w, v^{\comp }) \notin \mu ({\underline{\Delta }}')$. Write
$w = w_1+w_2$ where $w_i \in W_i$. Then by the argument in the preceding
paragraph we find that $Y^hw_1 = v^{\comp }$. Since $W_2$ is
$\gamma $-orthogonal to $W_1$ and $v^{\comp } \in W_1$, we obtain
\begin{displaymath}
\gamma (w, v^{\comp }) = \gamma (w_1,v^{\comp}) + \gamma (w_2, v^{\comp })
\, = \, \gamma (w_1, v^{\comp }).
\end{displaymath}
But $\gamma (w_1, v^{\comp }) \in \mu ({\underline{\Delta }}')$ by definition.
This is a contradiction. Hence $\mu ({\underline{\Delta }}') =
\mu (\underline{\Delta })$. \hfill $\square $ 

\begin{lem}\label{Lemma 1.3} 
We may write
$\underline{\Delta } = {\underline{\Delta }}' + \Delta $ where the
distinguished type ${\underline{\Delta }}'$ is indecomposable and 
nilpotent.
\end{lem} 

\noindent \textbf{Proof.}
If the distinguished type ${\underline{\Delta }}'$ is not
indecomposable, we find another distinguished type ${\underline{\Delta }}''$
of the same distinguished height and parameters and a type ${\Delta }'$
such that ${\underline{\Delta }}' = {\underline{\Delta }}'' +
{\Delta }'$, where $\dim {\Delta }' > 0$. Because 
$\dim {\underline{\Delta }}'
> \dim {\underline{\Delta }}''$ after a finite number of repetitions, we
obtain a distinguished type ${\underline{\widetilde{\Delta }}}$
which we can no longer write as a sum of a distinguished type and a type,
namely,  $\underline{\Delta } = {\underline{\widetilde{\Delta }}} +
\widetilde{\Delta }$. In other words,
$\underline{\widetilde{\Delta }}$ is an indecomposable
distinguished type. By lemma \ref{Lemma 1.2}
it has the same distinguished height and
parameters as the distinguished type $\underline{\Delta }$.
\par
We now show that the indecomposable distinguished type
$\underline{\widetilde{\Delta }}$, represented by $(Y|W, W, v^{\comp })$, is
nilpotent. Let $W_0$ be the generalized eigenspace of $Y|W$ corresponding
to the eigenvalue $0$. Then $W_0$ is $Y$-invariant,
$\gamma $-nondegenerate and contains
$v^{\comp }$. On $W_0$ the linear map $Y$ is nilpotent.
{}From the fact that the distinguished type $\underline{\widetilde{\Delta }}$
is indecomposable, it follows that the triple
$(Y|W_0, W_0, v^{\comp }; \gamma )$ equals the triple
$(Y|W, W, v^{\comp }; \gamma )$. Hence the indecomposable distinguished type
$\underline{\widetilde{\Delta }}$ is nilpotent. \hfill $\square $

\section{Indecomposable distinguished types}\label{section 3}

In this section we classify indecomposable distinguished types.
We start by giving a rough description of the possible indecomposable
distinguished types, which we then refine to a classification. \medskip

\noindent Let $\underline{\Delta }$ be a distinguished type. There are two
cases:
\par
1. the set of parameters $\mu (\underline{\Delta })$ contains
a nonzero parameter;
\par \noindent
or
\par
2. $\mu (\underline{\Delta }) = \{ 0 \} $. \medskip

\noindent \textsc{Case} 1. Suppose that the triple $(Y, W, v^{\comp })$
represents the distinguished type $\underline{\Delta }$, which we
assume has distinguished height $h$. Using lemma \ref{Lemma 1.3} write
$\underline{\Delta } = {\underline{\Delta }}' + \Delta $, where
${\underline{\Delta }}'$ is an indecomposable distinguished type of
distinguished height $h$ represented by $(Y|W_1, W_1, v^{\comp })$ 
with $W_1$ a
$\gamma $-nondegenerate, $Y$-invariant subspace of $W$ which contains
$v^{\comp }$.
Choose $w\in W_1$ so that $Y^hw= v^{\comp }$ and $\gamma (w,v^{\comp }) =
\mu  \ne 0$.\footnote{This implies that $h$ is even. Suppose not. Then
\begin{displaymath}
\gamma (w, Y^hw) = (-1)^h\, \gamma (Y^hw,w) \, =  -
\gamma (w, Y^hw),
\end{displaymath}
since $\gamma $ is symmetric. Hence $\gamma (w, Y^hw) =0$, which
is a contradiction.}
Look at the subspace
\begin{displaymath}
\widetilde{W} = \spann \{ w,\, Yw, \, \ldots \, , Y^hw \}
\end{displaymath}
 of $W$. Clearly $v^{\comp } \in \widetilde{W}$. On $\widetilde{W}$ consider
the $(h+1)\times (h+1)$ Gram matrix $G = \big( \gamma (Y^iw, Y^jw) \big) =
\big(\pm \gamma (w, Y^{i+j}w) \big)$, since 
$Y \in {\oo (W, \gamma )}_{v^{\comp }}$.
Because $Y^{h+1}w=Yv^{\comp }=0$, we have
$Y^{h+1}|\widetilde{W} =0$. Therefore, all the entries of $G$ below the
antidiagonal are $0$. On the other hand, because
\begin{displaymath}
\gamma (Y^iw, Y^{h-i}w) = \pm \gamma (w, Y^hw) \, = \, \pm \mu \ne 0,
\end{displaymath}
all the entries of $G$ on the antidiagonal are nonzero. Hence
$\det G \ne 0$, that is, $\widetilde{W}$ is $\gamma $-nondegenerate.
As ${\underline{\Delta }}'$ was assumed to be indecomposable, it follows
that $W_1 = \widetilde{W}$. Note that $(Y|\widetilde{W}, \widetilde{W},
v^{\comp })$ has one
Jordan block and therefore ${\underline{\Delta }}'$ is uniform.
This completes case 1. \medskip

\noindent \textsc{Case} 2. Suppose that the triple $(Y,W,v^{\comp })$
represents the distinguished type $\underline{\Delta }$, which we
assume has distinguished height $h$. Using lemma \ref{Lemma 1.3} write
$\underline{\Delta } = {\underline{\Delta }}' + \Delta $, where
${\underline{\Delta }}'$ is a nilpotent indecomposable distinguished type of
distinguished height $h$ represented by $(Y|W_1, W_1, v^{\comp })$ with $W_1$ a
$\gamma $-nondegenerate, $Y$-invariant subspace of $W$ which
contains $v^{\comp }$. Consider the pair $(Y|W_1, W_1)$
and  the type $\widetilde{\Delta }$ which it
represents. From the results of
\cite{burgoyne-cushman} we may write $\widetilde{\Delta } =
{\Delta }_{1} + \cdots +\, {\Delta }_{r}$, where ${\Delta }_{j}$
are indecomposable types uniform of height $h_j$, sorted so that
$h_1 \le h_2 \le \cdots \le h_r$.
Suppose that $(Y|W_{j}, W_{j})$ represents ${\Delta }_{j}$.
Then $v^{\comp }$ is a sum of its components in the $W_{j}$, but some of
those components may be zero. Let
$\widehat{W} = W_{k}$ where $k$ is the smallest index such that
$v^{\comp }$ has a nonzero component $\widehat{v^{\comp }}$ in $\widehat{W}$.
Consider the type $(Y|\widehat{W}, \widehat{W})$.
Then $Y|\widehat{W}$ annihilates $\widehat{v^{\comp }}$ and
the height of $(Y|\widehat{W}, \widehat{W})$ equals the
distinguished height $h$ of
${\underline{\Delta }}'$.
Choose $z \in \widehat{W}$ such that
$\gamma (z, v^{\comp }) =\gamma (z, \widehat{v^{\comp }}\,)\ne 0$.
This is possible since
$\widehat{W}$ is $\gamma $-nondegenerate.
Choose $w\in W_1$ so that $Y^hw= v^{\comp }$.
Consider the $Y$-invariant subspace
\begin{displaymath}
\widetilde{W} = \spann \{ w, \, Yw, \, \ldots \, Y^hw; \, z, \, Yz, \ldots
, \, Y^hz \} .
\end{displaymath}
Let $n = h+1$. Note that
$Y^{h+1}|\widetilde{W} =0$ and $\gamma (z, Y^hw) \ne 0$ by definition of $z$
and $w$. Moreover $\gamma (w, Y^hw) =0$ since $\mu (\underline{\Delta }) =
\{ 0 \} $ by hypothesis. Look at the $2n \times 2n$ Gram matrix
\begin{displaymath}
G = \mbox{{\footnotesize $ \left( \begin{array}{c|c}
g_{i,j} & g_{i,j+n} \\ \hline
g_{i+n,j} & g_{i+n,j+n} \end{array} \right) $}} =
\mbox{{\footnotesize $ \left( \begin{array}{c|c}
\gamma (Y^{i-1}w,Y^{j-1}w ) & \gamma (Y^{i-1}w,Y^{j-1}z) \\ \hline
\gamma (Y^{i-1}z, Y^{j-1}w ) & \gamma (Y^{i-1}z, Y^{j-1}z)
\end{array} \right) .$}}
\end{displaymath}
The entries of $G$ satisfy the following conditions: i)
$g_{i,j} = g_{n+i,j} \, = \, g_{i,n+j} \, = \, g_{n+i,n+j} \, = \, 0$,
when $i+j \ge n+2$ and $1 \le i,j \le n$; ii) $g_{i,j+n} = 
g_{i+n, j} \ne 0$, where $i+j = n+1$; iii) $g_{i,j} =0$, where $i+j = n+1$. 
Thus $G$ has its nonzero entries on or above the antidiagonal
of each $n \times n$ block except the upper left hand one, 
where even the antidiagonal elements are zero. Thus the matrix 
$G$ has the form
\begin{displaymath}
\mbox{{\footnotesize $\left( \begin{array}{ccc|ccc}
\mbox{{\Large $\ast $}}   & & 0      &  \mbox{{\Large $\ast $}}   & & + \\
& &  &   & &   \\
0      & & \mbox{{\Large $0$}}       & +      & & \mbox{{\Large $0$}} \\
\hline
\raisebox{-2pt}{\mbox{{\Large $\ast $}}}   & & + &
\raisebox{-2pt}{\mbox{{\Large $\ast$}}}   & & \ast \\
  & &  &  & &      \\
+      & & \mbox{{\Large $0$}}  & \ast   & & \mbox{{\Large $0$}}
\end{array} \right) $, }}
\end{displaymath}
where $+$ denotes a nonzero entry. Expanding $\det G$ by minors of the
${h+1}^{\mathrm{st}}$ column, one sees that $\det G$ is
a nonzero number times the $[h+2, h+1]$ minor. Expanding this 
minor by its last
column gives a nonzero number times
a matrix with the same form as the original $G$ but with one fewer
row and column. Clearly when $G$ is a $2\times 2$, we have $\det G \ne 0$.
By induction we have 

\begin{lem} \label{Lemma 2.1} $\det G = \pm \prod^{2n}_{k=1} g_{k,2n-k+1} 
\ne 0$. 
\end{lem} 

Thus $\widetilde{W}$ is a $2h+2$-dimensional, $Y$-invariant,
$\gamma $ nondegenerate
subspace of $W_1$, which contains the vector $v^{\comp }$. Since
${\underline{\Delta }}' $ is indecomposable, the triple
$(Y|\widetilde{W}, $ $ \widetilde{W}, v^{\comp })$ represents the
distinguished type ${\underline{\Delta }}' $. Note that
${\underline{\Delta }}' $ is made up of two Jordan blocks of size
$h+1$ and hence is uniform. This completes case 2 of
the rough description of indecomposable distinguished types.
\hfill $\square $ \medskip

We now classify indecomposable distinguished types. 

\begin{prop} \label{Proposition 2.2} 
Let $\underline{\Delta }$ be
an indecomposable distinguished type of distinguished height $h$, which
is represented by the triple $(Y,W, v^{\comp })$. Then exactly one of
the following alternatives holds.
\begin{itemize}
\item[1.] $h$ is \emph{even}, $h>0$, and there is a basis
\begin{equation}
\mbox{{\tiny $ \{ w, \, Yw, \ldots \, , Y^{h/2-1}w; \,
\varepsilon \, Y^hw, -\varepsilon \, Y^{h-1}w, \ldots \, ,
(-1)^{h/2-1}\varepsilon Y^{h/2+1}w; \, Y^{h/2}w \} $}} ,
\label{eq-s2four}
\end{equation}
where the Gram matrix of $\gamma $ is {\tiny $ \begin{pmatrix}
0 & I_{h/2} & 0 \\
I_{h/2} & 0 & 0 \\
0 & 0 & (-1)^{h/2}\varepsilon \end{pmatrix} $} and
$v^{\comp } = \mu \, Y^h w$ with $\mu >0$. We call $\mu $ a \emph{modulus}.
Here
${\varepsilon }^2 =1$. We use the
notation ${\underline{\Delta }}^{\varepsilon }_{h}(0), \, \, \mu $.
\item[2.] $h$ is \emph{odd} and there is a
basis
\begin{equation}
\{ Y^hz, \, -Y^{h-1}z, \ldots \, , (-1)^hz\, ;w, \, Yw, \,
 \ldots \, , Y^hw \} ,
\label{eq-s2five}
\end{equation}
where the Gram matrix of $\gamma $ is {\tiny $\begin{pmatrix}
0 & I_{h+1} \\
I_{h+1} & 0
\end{pmatrix} $} and
$v^{\comp } =  Y^hw$. We use the
notation ${\underline{\Delta }}_h(0,0)$.
\item[3.] $h$ is \emph{even} and there is a
basis
\begin{equation}
\{ Y^hz, \, -Y^{h-1}z, \ldots \, , (-1)^hz\, ;w, \, Yw, \,
 \ldots \, , Y^hw \} ,
\label{eq-s2six}
\end{equation}
where the Gram matrix of $\gamma $ is {\tiny $\begin{pmatrix}
0 & I_{h+1} \\
I_{h+1} & 0
\end{pmatrix} $} and
$v^{\comp } =  Y^hw$. We use the
notation $\underline{{\Delta }^{+}_h(0)+{\Delta }^{-}_h(0)}$ .
\end{itemize}
\end{prop}

\noindent \textbf{Proof.} Using our rough classification of distinguished
indecomposable types, let us prove the proposition.  \medskip

Suppose that we are in case 1 of the rough classification. Then
$\underline{\Delta }$ is represented by
the triple $(Y,W, v^{\comp})$ where $W = \spann \{ w, \, Yw,
\, \ldots \, , Y^hw \} $
and $\gamma (w, Y^hw) \ne 0$. Hence $h$ is even and $h>0$ because
$v^{\comp}$ is isotropic, while $\gamma (w, Y^hw) \ne 0$.
Since $\underline{\Delta }$
is uniform we may form $\overline{W} = W/YW$. Clearly, $\dim \overline{W} =1$.
On $\overline{W}$ the
inner product $\gamma $ induces a symmetric bilinear form $\overline{\gamma }$
defined by $\overline{\gamma }(\overline{v}, {\overline{v}}') =
\gamma (v, Y^hv')$. Since $\gamma (w, Y^hw) \ne 0$,
the vector $\overline{w}$ is nonzero
and forms a basis of $\overline{W}$. Rescaling, we may assume that
$\overline{\gamma }(\overline{w}, \overline{w}) = \varepsilon $,
where ${\varepsilon }^2 =1$. By \cite[prop.\ 2, p.343]{burgoyne-cushman}
any uniform type is determined by its height and its $(\overline{W},
\overline{\gamma })$, so
we may choose a vector $w \in W$ which generates the
basis
(\ref{eq-s2four}) of case 1 of the proposition, $\gamma $-adapted in the sense
that its Gram matrix is as indicated in the proposition.
Indeed such a $\gamma $-adapted basis describes a type that has
the required height and  $(\overline{W},
\overline{\gamma })$. In terms of this basis
there is a unique nonzero number $\mu $ such that $v^{\comp } =
\mu \,Y^h w$. Replacing $w$ with $-w$, if necessary, we can assume
that $\mu >0$. We call $\mu $ a modulus.
We compute that
\begin{displaymath}
{\gamma }(\mu w,{{v}}^{\comp }) =
\overline{\gamma }(\mu \overline{w},\mu  \overline{w}) = {\mu }^2 \,
\varepsilon ,
\end{displaymath}
which shows that $\mu (\underline{\Delta }) = \{ {\mu }^2\varepsilon \}$.
Thus $\mu (\underline{\Delta })$ determines $\mu $ and $\varepsilon $.
So $\underline{\Delta }$ is
a distinguished indecomposable type made up of one Jordan block. Moreover,
we have $\dim \underline{\Delta } =h+1$, $\mathrm{ind}\, \underline{\Delta }
= ${\tiny $\left\{ \begin{array}{rl}
h/2, & \mbox{if $(-1)^{h/2}\varepsilon =1$} \\
h/2+1, & \mbox{if $(-1)^{h/2}\varepsilon =-1$} \end{array} \right. $} and
$\underline{\Delta }$ has distinguished height $h$
and a unique modulus $\mu >0$.
The type of $(Y,W)$ is denoted ${{\Delta }}^{\varepsilon }_{h}(0)$ in
\cite{burgoyne-cushman}.
\medskip

Now suppose that we are in case $2$ of the rough classification.
Then the distinguished type
$\underline{\Delta }$ of distinguished height $h$ is represented by the
triple $(Y,W, v^{\comp })$ with
\begin{displaymath}
W = \spann \{ w,\, Yw, \ldots \, , Y^hw, \, z, \, Yz,\, \ldots \, , Y^hz \} ,
\end{displaymath}
and $v^{\comp } = Y^hw$. Moreover, $\gamma (w,v^{\comp }) =0$ and
$\gamma (z, v^{\comp }) \ne 0$. There are two subcases. \medskip

\noindent Suppose that
$h$ is odd. Since $\underline{\Delta }$ is uniform, we may form
$\overline{W} = W/YW$. On $\overline{W}$ the inner product $\gamma $
induces a skew symmetric bilinear form $\overline{\gamma }$ defined by
$\overline{\gamma }(\overline{v}, {\overline{v}}') = \gamma (v, Y^hv')$.
Clearly, $\overline{W} = \spann \{ \overline{w}, \overline{z} \} $ and
from $\overline{\gamma }(\overline{w}, {\overline{z}})\ne 0$ it follows that
$\overline{W}$ is $\overline{\gamma }$ nondegenerate.
Up to isomorphism there is only one nondegenerate
skew symmetric bilinear form of dimension two,
and it is indecomposable.
So $\overline{W}$ is $\overline{\gamma }$ indecomposable.
Using \cite[prop.\ 2, p.343]{burgoyne-cushman}
again we may choose vectors $w,z \in W$ which generate the
$\gamma $-adapted basis
(\ref{eq-s2five}) of case 2 of the proposition. We now need to
show that we can choose this basis so that $v^{\comp } = Y^hw$.
We know that $v^{\comp } = \alpha \, Y^hw + \beta \, Y^hz$ is a nonzero
vector in $\ker Y|W$. If $\alpha \ne 0$, let 
{\footnotesize $\begin{pmatrix} w' \\ z' \end{pmatrix}$}$ = ${\footnotesize $\begin{pmatrix} \alpha & \beta \\
0 & 1/\alpha \end{pmatrix} $}$\, ${\footnotesize $\begin{pmatrix} w \\ z 
\end{pmatrix} $}; while if $\alpha =0$ let 
{\footnotesize $\begin{pmatrix} w' \\ z' \end{pmatrix} $}$ = ${\footnotesize $\begin{pmatrix} 0 & \beta  \\
-1/\beta & 0 \end{pmatrix} $}$ \, ${\footnotesize $\begin{pmatrix} w \\ z 
\end{pmatrix} $}. We rewrite the definition as 
{\tiny  $\begin{pmatrix} w' \\ z' \end{pmatrix}  = 
\begin{pmatrix} a & b \\ c & d\end{pmatrix}  \, \begin{pmatrix} w \\ z 
\end{pmatrix}$,} where $ad-bc=1$. We now show that $w'$ and $z'$ generate
the $\gamma $-adapted basis $\{ Y^hz', \, -Y^{h-1}z', \ldots \, , (-1)^hz'\, ;w', \, Yw', \, \ldots \, , Y^hw' \} $ of $W$. This follows because for every
$j$ between $0$ and $h$ we have
\begin{displaymath}
\gamma (Y^iw', Y^jw') = \gamma (Y^iz', Y^jz') \, = \, 0
\end{displaymath}
and
\begin{eqnarray*}
\gamma (Y^jw', (-1)^jY^{h-j}z') & = &
(-1)^j\gamma (Y^j(aw+bz), Y^{h-j}(cw+dz)) \\
& = &\gamma (aw+bz,Y^h(cw+dz)) \\
& = & ac \, \overline{\gamma } (\overline{w}, \overline{w}) +
bd\, \overline{\gamma } (\overline{z}, \overline{z}) +
(ad - bc)\overline{\gamma } (\overline{w}, \overline{z}) \\
& = & \overline{\gamma } (\overline{w}, \overline{z})  \, = \, 1.
\end{eqnarray*}
By construction
$v^{\comp } = Y^h w'$.
Summarizing, we have shown that
$\underline{\Delta }$ is
a distinguished indecomposable type made up of two Jordan blocks. Also
$\dim \underline{\Delta } =2(h+1)$, $\mathrm{ind}\, \underline{\Delta }
= h+1$ and $\underline{\Delta }$ has distinguished height $h$, which is odd.
The type of $(Y,W)$ is denoted ${{\Delta }}_{h}(0,0)$ in
\cite{burgoyne-cushman}.\medskip

\noindent Suppose that $h$ is even. Since $\underline{\Delta }$ is uniform,
we may form $\overline{W} = W/YW$. On $\overline{W}$ the inner product
$\gamma $ induces a symmetric bilinear form $\overline{\gamma }$ defined
by $\overline{\gamma }(\overline{v}, {\overline{v}}') = \gamma (v, Y^hv')$.
Since $\gamma (z, Y^hw) \ne 0$ by hypothesis, we see that
$\gamma (\overline{z}, \overline{w}) \ne 0$ and $\overline{W} =
\spann \{ \overline{z}, \overline{w} \} $. Therefore the
reduced type $(\overline{Y}, \overline{W}; \overline{\gamma })$ is
\emph{not} indecomposable. Since $\gamma (w, Y^hw) =0$, the vector
$\overline{w}$ is a nonzero and $\overline{\gamma }$-isotropic.
Let $\overline{y} = \frac{1}{\gamma (\overline{z}, \overline{w})}
\left( \overline{z} -
\frac{\gamma (\overline{z},\overline{z})}{2\gamma (\overline{z},\overline{w})}
\overline{w} \right) $. Then $\overline{y}$ is a
$\overline{\gamma }$-isotropic vector in $\overline{W}$ and
$\overline{\gamma }(\overline{y}, \overline{w}) = 1$. Thus the
matrix of $\overline{\gamma }$ with respect to the basis
$\{ \overline{y}, \overline{w} \} $ is {\tiny $
\begin{pmatrix} 0 & 1 \\ 1 & 0 \end{pmatrix} $}.
Using  \cite[prop.\ 2, p.343]{burgoyne-cushman} we may choose vectors
$\widetilde{w},\widetilde{z} \in W$ which generate the
$\gamma $-adapted basis
(\ref{eq-s2six}) of case 3 of the proposition. We now need to
show that we can choose this basis so that $v^{\comp } = Y^h\widetilde{w}$.
Since $v^{\comp } \in \ker Y|W$, we see that $v^{\comp } \in
\spann \{ Y^h \widetilde{w}, Y^h\widetilde{z} \} $. Now write
$v^{\comp } = Y^h (\alpha \widetilde{w} + \beta \widetilde{z}) $. As
$\gamma(\alpha
\widetilde{w} + \beta \widetilde{z},v^{\comp })=
2\alpha\beta\in \mu(\underline{\Delta }) \, = \{ 0 \}$,
we must have $\alpha=0$ or $\beta=0$.
 If $v^{\comp } = \alpha \, Y^h \widetilde{w}$, where
$\alpha \ne 0$, then put $z'={\alpha }^{-1}\widetilde{z}$,
$w'={\alpha }\widetilde{w}$.
If
$v^{\comp }= \beta \, Y^h \widetilde{z}$ with $\beta \ne 0$ then
put $z'={\beta }^{-1}\widetilde{w}$,
$w'={\beta }\widetilde{z}$.
In either case $v^{\comp } = Y^h{w'}$ and
\begin{displaymath}
\{ Y^hz', \, -Y^{h-1}z', \ldots \, , (-1)^hz'\, ;w', \, Yw', \,
 \ldots \, , Y^hw' \}
\end{displaymath}
is a basis of $W$ with respect to which
the matrix of $\gamma $ is {\tiny $\begin{pmatrix}
0 & I_{h+1} \\ I_{h+1} & 0 \end{pmatrix} $}. Note $\underline{\Delta }$
is a distinguished indecomposable type made up of two Jordan blocks.
Also $\dim \underline{\Delta } =2(h+1)$ with
$\mathrm{ind}\, \underline{\Delta }
= h+1$ and $\underline{\Delta }$ has distinguished height $h$, which is even.
The type of $(Y,W)$ is decomposable and
is denoted ${\Delta }^{+}_h(0)+{\Delta }^{-}_h(0)$ in
\cite{burgoyne-cushman}.

One may look at the above computation as exploiting the fact that there
is an action
of $\Oo(\overline{W},\overline{\gamma })$ on $\ker Y|W$. In the last
two cases the action has only one orbit
of nonzero isotropic vectors, while in the first case there are moduli.
The action can be understood in terms of the Jacobson Morozov theorem.
\medskip

The three cases are obviously exclusive. Note that one can
distinguish them
by $\mathrm{dht}(\underline{\Delta }) $ and $\mu (\underline{\Delta }) $.
This proves proposition~\ref{Proposition 2.2}.\,\hfill $\square $ \medskip

\noindent \textbf{Proof of theorem 1}
Let $\underline{\Delta }$ be a distinguished
type. By lemma \ref{Lemma 1.3} we may write
$\underline{\Delta } = \widetilde{\underline{\Delta }} + \Delta $ where the
distinguished
type $\widetilde{\underline{\Delta }}$ is indecomposable and nilpotent.
By the main result of \cite[theorem, p.343]{burgoyne-cushman} applied to
$\Delta $, we can write
\begin{equation}
\underline{\Delta } = \widetilde{\underline{\Delta }} +
{\Delta }_1 + \cdots + {\Delta }_r,
\label{eq-s2seven}
\end{equation}
where  ${\Delta }_i$ for $1 \le i \le r$ are
indecomposable types.
By lemma \ref{Lemma 1.2}  $\widetilde{\underline{\Delta }}$ is of the same
distinguished height and parameters as $\underline{\Delta }$.
 Suppose that $\underline{\Delta }$ has another
such decomposition, namely
\begin{equation}
\underline{\Delta } = {\widetilde{\underline{\Delta }}}' +
{\Delta }'_1 + \cdots + {\Delta }'_s,
\label{eq-s2eight}
\end{equation}
where ${\widetilde{\underline{\Delta }}}'$ is an indecomposable
distinguished type and ${\Delta }'_j$ for $1 \le j \le s$ are
indecomposable types. By lemma \ref{Lemma 1.2} the distinguished height,
say $h$, of
$\widetilde{\underline{\Delta }}$ and
${\widetilde{\underline{\Delta }}}'$ are the same. Say that
$\widetilde{\underline{\Delta }}$ and
${\widetilde{\underline{\Delta }}}'$ are represented by
the triples $(Y,W, v^{\comp })$ and $(Y', W', (v^{\comp })')$.
Suppose that $h$ is odd. Then the linear map $P:W \rightarrow W'$
for which $PY^iw = (Y')^iw'$ and $PY^iz = (Y')^iz'$ where
$0 \le i \le h$ and $w,z$ and $w',z'$ are
vectors given in the basis (\ref{eq-s2five}) of case 2 of proposition
\ref{Proposition 2.2}
is an equivalence
between the triples $(Y,W, v^{\comp })$ and $(Y', W', (v^{\comp })')$. Next
suppose that $h$ is even and that $(Y,W, v^{\comp })$ and
$(Y', W', (v^{\comp })')$ have one Jordan chain. Since by lemma \ref{Lemma 1.2}
the parameters of $\widetilde{\underline{\Delta }}$ and
${\widetilde{\underline{\Delta }}}'$ are the same, using the basis
(\ref{eq-s2four}) of case 1 of proposition \ref{Proposition 2.2}
we can again construct an
equivalence
between $\widetilde{\underline{\Delta }}$ and
${\widetilde{\underline{\Delta }}}'$. We can also handle the case
when $h$ is even and $(Y,W, v^{\comp })$ and
$(Y', W', (v^{\comp })')$ have two Jordan chains. Thus in
every case $(Y,W, v^{\comp })$ and
$(Y', W', (v^{\comp })')$ are equivalent, that is,
$\widetilde{\underline{\Delta }} = {\widetilde{\underline{\Delta }}}'$.
\medskip

Now we need only show that $r=s$ and ${\Delta }_i = {\Delta }'_i$. But this
follows from the the main result of \cite[theorem, p.343]{burgoyne-cushman},
because
 ${\Delta }_1 + \cdots {\Delta }_r$ and
${\Delta }'_1 + \cdots {\Delta }'_s$ are sums of indecomposable types,
while $\widetilde{\underline{\Delta }} = {\widetilde{\underline{\Delta }}}'$
implies that the underlying types of $\widetilde{\underline{\Delta }}$
and ${\widetilde{\underline{\Delta }}}'$ are equal.
This proves theorem 1. \hfill $\square $

\section{Adjoint orbits of the Poincar\'{e} group}\label{section 4}

In this section we use the above theory to determine the
orbits of the adjoint action of the \emph{Poincar\'{e}} group on its
Lie algebra. \medskip

Let $G = \mathrm{diag}(-1,-1,-1,1)$ be the matrix of
a Lorentz inner product on ${\R }^4$ with respect to the standard
basis $\{ e_1, \ldots \, , e_4 \} $. The Poincare group is
the affine Lorentz group, which is the semidirect product
$\Oo (3,1) \ltimes  {\R }^4$ of the Lorentz group
$\Oo (3,1) = \Oo ({\R }^4,G)$ with the abelian group
${\R }^4$. In \S \ref{section 1} we have shown that the Poincar\'{e} group
is the isotropy group ${\Oo ({\R }^6,K)}_{e_5}$ of the
orthogonal group $\Oo ({\R }^6, K)$, where the matrix of
the inner product $K$ with respect to the basis
$\{ e_0,e_1, \ldots \, , e_4,e_5 \} $ of ${\R }^6$ is standard. The Lie algebra of the Poincar\'{e} group is isomorphic to
the Lie algebra ${\oo ({\R }^6, K)}_{e_5}$ of ${\Oo ({\R }^6,K)}_{e_5}$. 
All the conjugacy classes in ${\oo ({\R }^6, K)}_{e_5}$
are given in table 3 below. \medskip

First we list all the possible ${\oo ({\R }^6, K)}_{e_5}$-indecomposable
distinguished types, meaning indecomposable distinguished types that may occur as summand of some 
$(Y,{\R }^6,e_5; K)$. \bigskip 

\noindent \hspace{.75in}\begin{tabular}{ccccc}
&\multicolumn{1}{c}{type (modulus $\alpha >0$)} &\multicolumn{1}{c}{dim} &
\multicolumn{1}{c}{index} & \multicolumn{1}{c}{$v^{\comp }$} \\
\hline
$1$. & $\underline{{\Delta }^{-}_4(0)}, \, \alpha >0$ & $5$ & $3$ & $\alpha
\, Y^4w$ \\
$2$. & $\underline{{\Delta }^{+}_4(0)}, \, \alpha >0$ & $5$ & $2$ & $\alpha
\, Y^4w$ \\
$3$. & $\underline{{\Delta }_1(0,0)}$ & $4$ & $2$ &  $Yw$ \\
$4$. & $\underline{{\Delta }^{+}_2(0)}, \, \alpha >0$ & $3$ & $2$ & $\alpha
\, Y^2w$ \\
$5$. & $\underline{{\Delta }^{-}_2(0)}, \, \alpha >0$ & $3$ & $1$ & $\alpha
\, Y^2w$ \\
$6$. & $\underline{{\Delta }^{+}_0(0)+{\Delta }^{-}_0(0)}$ & $2$ & $1$ &
$w$
\end{tabular}

\begin{center}Table 1. Possible ${\oo ({\R }^6, K)}_{e_5}$-indecomposable
distinguished types. \end{center} \bigskip  

Note we express $v^{\comp }$ using the basis given
in proposition \ref{Proposition 2.2}.  \medskip

We now show that all the possible indecomposable distinguished types are listed
in table 1. The possible eigenvalue combinations are $0\, 0$; $0$; and $0 + 0$.
Here, for instance, $0 + 0$ stands for a decomposable two
dimensional $(\overline{Y}, \overline{W}; \overline{\gamma })$ with eigenvalue
zero for each summand.
The corresponding heights and signs are $1$; $4^{\pm },\, 2^{\pm}$; and $0$.
So table 1 lists all the possibilities. \medskip

Next in table 2 below we list the possible 
${\oo ({\R }^6,K)}_{e_5}$-indecomposable types,
see \cite[table II, p.349]{burgoyne-cushman}. That is, we look for types that
occur as proper summand of some $(Y,{\R }^6; K)$. We do not claim they all 
actually occur in the setting of theorem 1. \bigskip 

\noindent \hspace{.35in}\begin{tabular}{cccc|cccc}
&\multicolumn{1}{c}{type} & \multicolumn{1}{c}{dim} &
\multicolumn{1}{c}{index} & &\multicolumn{1}{c}{type} &
\multicolumn{1}{c}{dim} & \multicolumn{1}{c}{index}  \\ \hline
$1$. & ${\Delta }^{-}_4(0) $ & $5$ & $3$
& $8$. & ${\Delta }^{-}_2(0) $ & $3$ & $1$  \\
$2$. & ${\Delta }^{+}_4(0) $ & $5$ & $2$
& $9$. & ${\Delta }^{-}_0(\zeta , \mathrm{IP}) $ & $2$ & $2$  \\
$3$. & ${\Delta }_0(\zeta, \mathrm{CQ}) $ & $4$ & $2$
& $10$. & ${\Delta }_0(\zeta, \mathrm{RP}) $ & $2$ & $1$   \\
$4$. & ${\Delta }_1(\zeta, \mathrm{RP}) $ & $4$ & $2$
& $11$. & ${\Delta }^{+}_0(\zeta , \mathrm{IP}) $ & $2$ & $0$   \\
$5$. & ${\Delta }^{\varepsilon}_1(\zeta , \mathrm{IP}) $ & $4$ & $2$
& $12$. & ${\Delta }^{-}_0(0) $ & $1$ & $1$   \\
$6$. & ${\Delta }_1(0,0) $ & $4$ &  $2$
& $13$. & ${\Delta }^{+}_0(0) $ & $1$ & $0$ \\
$7$. & ${\Delta }^{+}_2(0) $ & $3$ & $2$  &&&&
\end{tabular} \bigskip 

\noindent \hspace{.75in} Table 2. Possible ${\oo ({\R }^6, K)}_{e_5}$-indecomposable
types. \bigskip 

\noindent Note in table 2 we have used the notation
${\Delta }_m(\zeta, \mathrm{CQ})  =
{\Delta }_m(\zeta, -\zeta , \overline{\zeta }, -\overline{\zeta })$,
$\zeta \ne \pm \overline{\zeta }$,
${\Delta }_m(\zeta, \mathrm{RP})  = {\Delta }_m(\zeta, -\zeta)$,
$\zeta = \overline{\zeta } \ne 0$,
${\Delta }_m(\zeta, \mathrm{IP})  = {\Delta }_m(\zeta, -\zeta)$,
$\zeta = -\overline{\zeta } \ne 0$, where $\zeta $ is the complex
eigenvalue of $Y$ with $( Y,W;K)$ a
representative of the $\oo ({\R }^6, K)$-indecomposable type.
For instance, ${\Delta }_m(\zeta, -\zeta , \overline{\zeta }, -\overline{\zeta })$
has height $m$ and four eigenvalues on $\overline{W}$.
\medskip

We now show that all the possible
$\oo ({\R }^6, K)_{e_5}$-indecomposable types are
listed in table 2. For each eigenvalue combination we have
the following possibilities for the heights and the signs,
because the dimension is at most five.
\begin{displaymath}
\begin{array}{r|c|c|c|c|c}
\mbox{eigenvalues} & \mathrm{CQ}  & \mathrm{IP} &  \mathrm{RP} &
0 & 0\, 0 \\ \hline
\rowspace \mbox{height and sign} & 0 & 0^{\pm},\,
1^{\pm }& 0, \, 1 & 0^{\pm}, \,
2^{\pm }, \, 4^{\pm } & 1.
\end{array}
\end{displaymath}
This gives a total of fourteen cases, two of which are covered
by case 5. Thus table 2 is complete. \medskip

Next we combine a given distinguished type in table 1 with a sum
of $\oo ({\R }^6, K)_{e_5}$-indecomposable types from table 2
so that their dimensions add up to $6$ and their indices add up to $4$. This 
gives the entries in table 3. \bigskip

\noindent \hspace{-.1in}\begin{tabular}{cccccc}
&\multicolumn{1}{c}{indecomposable}  &
& &  &  \\
&\multicolumn{1}{c}{distinguished type} &
\multicolumn{1}{c}{sum of ${\oo ({\R }^6, K)}_{e_5}$} &  &  \\
&\multicolumn{1}{c}{(modulus $\alpha >0$)} &
\multicolumn{1}{c}{indecomposable types} &\multicolumn{1}{c}{dim} &
\multicolumn{1}{c}{index} \\ \hline
$1$. &  $\underline{{\Delta }^{-}_4(0)}, \, \alpha  $ & & $5$ & $3$  \\
\rule{10pt}{0pt}a. & & $+{\Delta }^{-}_0(0)$ & $1$ & $1$ \\ \hline
$2$. & $\underline{{\Delta }_1(0,0)} $ & & $4$ & $2$  \\
\rule{10pt}{0pt}a. & & $+{\Delta }^{-}_0(\zeta , \mathrm{IP})$ & $2$ & $2$ \\
\rule{10pt}{0pt}b. & & $+{\Delta }^{-}_0(0)+ {\Delta }^{-}_0(0)$ & $2$ & $2$ \\
\hline
$3$. &  $\underline{{\Delta }^{+}_2(0)}, \, \alpha  $  & & $3$ & $2$  \\
\rule{10pt}{0pt}a. & & $+{\Delta }^{+}_2(0)$ & $3$ & $2$ \\
\rule{10pt}{0pt}b. & & $+{\Delta }^{-}_0(\zeta , \mathrm{IP}) +
{\Delta }^{+}_0(0)$ &
$3$ & $2$ \\
\rule{10pt}{0pt}c. & & $+{\Delta }_0(\zeta , \mathrm{RP})
+ {\Delta }^{-}_0(0)$ &
$3$ & $2$ \\
\rule{10pt}{0pt}d. & & $+{\Delta }^{-}_0(0)+{\Delta }^{-}_0(0)+
{\Delta }^{+}_0(0)$ & $3$ & $2$ \\ \hline
$4$. &  $\underline{{\Delta }^{-}_2(0)}, \, \alpha $  & & $3$ & $1$  \\
\rule{10pt}{0pt}a. & & $+{\Delta }^{-}_0(\zeta , \mathrm{IP}) +
{\Delta }^{-}_0(0)$ &
$3$ & $3$ \\
\rule{10pt}{0pt}b. & & $+{\Delta }^{-}_0(0)+{\Delta }^{-}_0(0)+
{\Delta }^{-}_0(0)$ &
$3$ & $3$ \\ \hline
$5$. &   $\underline{{\Delta }^{+}_0(0)+{\Delta }^{-}_0(0)}$ & & $2$ & $1$ \\
\rule{10pt}{0pt}a. & &  $+{\Delta }^{+}_2(0) + {\Delta }^{-}_0(0)$ & $4$ & $3$
\\
\rule{10pt}{0pt}b. & & $+{\Delta }^{-}_0(\zeta , \mathrm{IP})+
{\Delta }_0(\zeta , \mathrm{RP}) $ & $4$ & $3$ \\
\rule{10pt}{0pt}c. & & $+{\Delta }^{-}_0(\zeta , \mathrm{IP})+
{\Delta }^{-}_0(0)+ {\Delta }^{+}_0(0)$ & $4$ & $3$ \\
\rule{10pt}{0pt}d. & & $+{\Delta }_0(\zeta , \mathrm{RP})+{\Delta }^{-}_0(0)+
{\Delta }^{-}_0(0)$ & $4$ & $3$ \\
\rule{10pt}{0pt}e. & & $+{\Delta }^{-}_0(0)+{\Delta }^{-}_0(0)+
{\Delta }^{-}_0(0) +{\Delta }^{+}_0(0)$
& $4$ & $3$ \\
\end{tabular} \bigskip 

\noindent \hspace{1in}Table 3. Conjugacy classes in ${\oo ({\R }^6,K)}_{e_5}$. \bigskip 

The following list of dimension-index pairs shows that all the
${\Oo ({\R }^6, K)}_{e_5}$-conjugacy classes in ${\oo ({\R }^6,K)}_{e_5}$
are given in table 3. \medskip

\begin{center}
\begin{tabular}{ccl}
& &  \multicolumn{1}{c}{dimension-index pairs in}  \\
&\multicolumn{1}{c}{dimension-index pair} &
\multicolumn{1}{c}{sum of indecomposable types}  \\ \hline
$1$. & $(5,3)$ & $(1,1)$ \\
$2$. & $(4,2)$ & $(2,2)$, \,  $(1,1)+(1,1)$ \\
$3$. & $(3,2)$ & $(3,2)$, \, $(2,2)+(1,0)$, \, $(2,1)+(1,1)$, \\
&  &  $(1,1)+(1,1)+(1,0)$ \\
$4$. & $(3,1)$ & $(2,2)+(1,1)$, \,  $(1,1)+(1,1)+(1,1)$ \\
$5$. & $(2,1)$ & $(3,2)+(1,1)$, \, $(2,2)+(2,1)$, \\
& & $(2,2)+(1,1)+(1,0)$, $(2,1)+(1,1)+(1,1)$, \\
& &  $(1,1)+(1,1)+(1,1)+(1,0)$.
\end{tabular}
\end{center}

Below we show how to find explicit normal forms from
the decomposition into an indecomposable distinguished
${\oo ({\R }^6,K)}_{e_5}$-type and a sum of indecomposable
${\oo ({\R }^6,K)}_{e_5}$-types given in table 3. We do this for one case just to give the idea. \medskip

\begin{example}   $\underline{{\Delta }^{-}_4(0)}, \, \alpha
+{\Delta }^{-}_0(0)$. 
\end{example}

\noindent Write ${\R }^6 = V_1 \oplus V_2$, where $V_1$ and $V_2$ are
$Y$-invariant, $K$-orthogonal, ${\oo ({\R }^6,K)}_{e_5}$-indecomposable
subspaces where $(V_1, Y|V_1) \in {\Delta }^{-}_4(0)$, $e_5\in V_1$,
and $(V_2,$ $ Y|V_2) \in
{\Delta }^{-}_0(0)$. Now $Y=N$ is nilpotent on $V_1$ and $V_2$. Choose
a basis
\begin{displaymath}
\{ v_1, \, Nv_1, \, -N^4v_1,\,  N^3v_1;\,  N^2v_1 \}
\end{displaymath}
of $V_1$ as in case 1 of proposition \ref{Proposition 2.2}.
 Note that
$v^{\comp }= \alpha N^4v_1$ with $\alpha >0$. Also there is
a vector $v_2$ in $V_2$ such that $K(v_2, v_2) =-1$. With
respect to the  basis
\begin{displaymath}
\{ e_0, \ldots , \, e_5 \} =
\{ -{\alpha }^{-1}v_1, \, \onehalf Nv_1 -  N^3v_1,
\,  N^2v_1, \, v_2, \, \,\onehalf Nv_1+ N^3v_1; \, \alpha  N^4v_1 \}
\end{displaymath}
the matrix of $K$ is standard while the matrix of
$Y\in {\oo ({\R }^6,K)}_{e_5}$ is
\begin{displaymath}
\mbox{{\footnotesize $\left( \begin{array}{c|crrc|c}
0                 & 0         &  0         &  0       & 0 & 0 \\ \hline
-{\alpha }^{-1}   & 0         &  -\onehalf & 0        & 0 & 0 \\
0                 & \onehalf  &  0         & 0  & \onehalf & 0 \\
0                 & 0         &  0         & 0       & 0 & 0 \\
-{\alpha }^{-1}          & 0         &  \onehalf         & 0 & 0 & 0  \\ \hline
0           & -{\alpha }^{-1}  &  0  & 0 & {\alpha }^{-1} & 0 \end{array}
\right) $},}
\end{displaymath}
which is the desired normal form.

\section{Classification of coadjoint orbits}\label{section 5}
Our next aim is to determine a representative of each
orbit of the coadjoint action
\begin{displaymath}
{\Oo ({\R }^6, K)}_{e_5} \times {\oo ({\R }^6,K)}^{\ast }_{e_5}
\rightarrow {\oo ({\R }^6,K)}^{\ast }_{e_5}: (P, Y^{\ast })
\mapsto Y^{\ast }\comp \Ad_{P^{-1}}
\end{displaymath}
of the Poincar\'{e} group ${\Oo ({\R }^6, K)}_{e_5}$ on the dual
${\oo ({\R }^6,K)}^{\ast }_{e_5}$ of its Lie algebra.
More generally, we classify the coadjoint orbits of an affine orthogonal group.
As before, it is essential to our method that the affine orthogonal group
is viewed as an isotropy subgroup. Instead of types we will now employ cotypes.

\medskip

As always, the pair $(V,\gamma )$ is a finite dimensional real vector space
with
a nondegenerate inner product $\gamma $.
When $K$ is the Gram matrix of $\gamma$ with respect to some basis, we
often write $K$ for $\gamma$. For a vector
$v$ in $V$ let $v^{\ast }$ be the
linear function on $V$ given by $w \mapsto \gamma (v, w)$.
A \emph{tuple} $(V,Y,v;\gamma )$ is a pair $(V,\gamma )$,
a real linear map $Y \in \oo (V, \gamma )$ and a vector $v \in V$.
On the collection of all tuples we say that the tuple $(V, Y,v;\gamma )$ is
\emph{equivalent} to the tuple $(V', Y',v';\gamma ' )$ if and
only if there is a bijective real linear map $P:V \rightarrow V'$
such that (\textit{i}) $P^{\ast }{\gamma }' = \gamma $,
(\textit{ii}) $Pv =v'$, and
(\textit{iii}) there is a vector $w \in V$ such that
$Y' = P(Y + L_{w,v})P^{-1}$, where $L_{w,v} = w \otimes v^{\ast }
- v \otimes w^{\ast }$. 

\begin{fact}\label{Fact 1.} 
$L_{w,v} \in \oo (V, \gamma )$. \hfill $\Box$ 
\end{fact} 

\begin{fact} \label{Fact 2.} If $P \in \Oo (V, \gamma )$, then
$PL_{w,v}P^{-1} = L_{Pw,Pv}$.\hfill $\Box $ 
\end{fact} 

Being equivalent is an equivalence relation on the collection of
tuples. An equivalence class is a \emph{cotype}, which is denoted
by $\nabla $. If $(V,Y, v; \gamma )$ is a representative of $\nabla $, then
define the \emph{dimension} of $\nabla $ to be $\dim V$ and denote it by
$\dim \nabla $. Clearly, the notion of dimension is well defined.
A cotype is \emph{affine} if it has a representative
$(V, Y, v; \gamma )$, where $v$ is a nonzero, $\gamma $-isotropic
vector. \medskip

Suppose that we are in the situation of \S \ref{section 1}, where
$V = \R \times \widetilde{V} \times \R $ is a real vector space
with nondegenerate inner product $\gamma $ defined by
\begin{displaymath}
\gamma ((x, \widetilde{v},y),(x', {\widetilde{v}}' ,y')) =
\widetilde{\gamma }(\widetilde{v}, {\widetilde{v}}') +x'y+y'x,
\end{displaymath}
where $\widetilde{\gamma }$ is a nondegenerate inner product on
$\widetilde{V}$.
Suppose that with respect to the standard basis $\mathfrak{e} =
\{ e_0, e_1, \ldots , e_n, e_{n+1} \} $ of $V$ the matrix of
$\gamma $ is $K=${\tiny $\begin{pmatrix} 0 & 0 & 1 \\ 0 & G & 0 \\
1 & 0 & 0 \end{pmatrix} $}, where $G$ is the matrix of
$\widetilde{\gamma }$ with respect to the basis $\widetilde{\mathfrak{e}}=
\{ e_1, \ldots e_n \} $ of $\widetilde{V}$. \medskip

The following proposition explains the relevance of affine cotypes.
See also proposition \ref{prop5.4} below.  \medskip

For $Y \in \oo (V, \gamma )$ let ${\ell }_Y$ 
be the linear function on $\oo (V, \gamma )$ which
maps $Z$ to $\mathrm{tr}\, YZ$. Observe that the map 
$\oo (V,\gamma ) \rightarrow {\oo (V, \gamma)}^{\ast }:Y \mapsto {\ell }_Y$
is bijective. 

\begin{prop} \label{Prop5.1} 
The map 
\begin{equation}
(V,Y, e_{n+1}; K) \mapsto {\ell }_Y|{\oo (V, K)}_{e_{n+1}}
\label{eq-newnine}
\end{equation}
induces a bijection between affine cotypes on $(V,K)$ and
coadjoint orbits of ${\Oo (V,K)}_{e_{n+1}}$ on the dual
${\oo (V,K)}^{\ast }_{e_{n+1}}$ of its Lie algebra. 
\end{prop}  

\noindent \textbf{Proof.} The argument is a series of observations. \medskip

Suppose that the tuples $(V,Y, e_{n+1};K)$ and $(V,Y', e_{n+1};K)$ are
equivalent. Then there is a real linear map $P\in {\Oo (V, K)}_{e_{n+1}}$
and a vector $w\in V$ such that $Y' = P(Y+L_{w,e_{n+1}})P^{-1}$. \medskip

\noindent \textbf{Observation 1.} The matrix of $L_{w,e_{n+1}}$ with
respect to the standard basis $\mathfrak{e}$ of $(V,K)$ is
{\tiny $\begin{pmatrix} w_0 & 0 & 0 \\
\widetilde{w} & 0 & 0 \\ 0 &-{\widetilde{w}}^TG & -w_0 \end{pmatrix} ,$} where $w = w_0e_0 + \widetilde{w} +w_{n+1}e_{n+1} \in V$. \medskip

\noindent \textbf{Proof.} We compute
\begin{eqnarray*}
L_{w, e_{n+1}}(e_0) & = & (e^T_{n+1}Ke_0)w -(w^TKe_0)e_{n+1} =
w-w_{n+1}e_{n+1} =\\& = & w_0 e_0 + \widetilde{w}; \\
L_{w, e_{n+1}}(e_i) & = & (e^T_{n+1}Ke_i)w -(w^TKe_i)e_{n+1} =
-({\widetilde{w}}^TGe_i)e_{n+1}; \\
L_{w, e_{n+1}}(e_{n+1}) & = & (e^T_{n+1}Ke_{n+1})w -(w^TKe_{n+1})e_{n+1} =
-w_0 e_{n+1}. \qquad \Box
\end{eqnarray*}

\noindent \textbf{Observation 2.} For $P \in \Oo (V,K)$ and
$Y \in \oo (V,K)$ we have
\begin{displaymath}
{\ell }_{PYP^{-1}} = {\Ad }^T_{P^{-1}}{\ell }_Y:={\ell }_Y\comp
{\Ad }_{P^{-1}}.
\end{displaymath}

\noindent \textbf{Proof.} Let $Z \in \oo (V, K)$. Then
\begin{eqnarray*}
{\ell }_{PYP^{-1}}(Z) & = & \mathrm{tr}\, \big( P(YP^{-1}Z) \big)
\, = \, \mathrm{tr}\, \big( (YP^{-1}Z)P  \big) \\
& = & \mathrm{tr} \, \big( Y(P^{-1}ZP) \big) \, = \,  {\ell }_Y({\Ad
}_{P^{-1}}Z) \\
& = & ({\Ad }^T_{P^{-1}}{\ell }_Y)Z. \qquad \Box
\end{eqnarray*}

\noindent \textbf{Observation 3.}
Let ${\oo (V,K)}^0_{e_{n+1}}$ be the set of all ${\ell }_X \in 
{\oo (V,K)}^{\ast }$ such that ${\ell }_X(Y) = \mathrm{tr}\, XY = 0$ 
for every $Y \in {\oo (V,K)}_{e_{n+1}}$. Then 
\begin{equation}
{\oo (V,K)}^0_{e_{n+1}} =\{ {\ell }_{L_{v,e_{n+1}}} \in 
{\oo (V,K)}^{\ast } \setrule \, v = v_0e_0 +\widetilde{v} \in V \} . 
\label{eq-newstar}
\end{equation}

\noindent \textbf{Proof.} With respect to the standard basis 
$\mathfrak{e}$ of $(V,K)$ the matrix of $X$ and $Y$ is 
\begin{displaymath}
\mbox{{\footnotesize $\left( \begin{array}{ccr}
x_0 & -{\widetilde{u}}^TG & 0 \\
\widetilde{x} & \widetilde{X} & \widetilde{u} \\
0 & -{\widetilde{x}}^TG & -x_0 \end{array} \right) $}} \, \, \, 
\mathrm{and} \, \, \, 
\mbox{{\footnotesize $\left( \begin{array}{ccc} 
0 & 0 & 0 \\ \widetilde{y} & \widetilde{Y} & 0 \\ 
0 & -{\widetilde{y}}^TG & 0 \end{array} \right) $}}, 
\end{displaymath}
respectively, where $x_0 \in \R$, $\widetilde{x}, \widetilde{y}, 
\widetilde{u} \in 
\widetilde{V}$, and $\widetilde{X}, \widetilde{Y} \in \oo (\widetilde{V},G)$. 
Suppose that ${\ell }_X \in  {\oo (V,K)}^0_{e_{n+1}}$. Then for every 
$Y \in {\oo (V,K)}_{e_{n+1}}$ 
\begin{align}
0 & = \mathrm{tr}(XY) = \mathrm{tr}\left[ \mbox{{\footnotesize $\left( 
\begin{array}{ccr}
x_0 & -{\widetilde{u}}^TG & 0 \\
\widetilde{x} & \widetilde{X} & \widetilde{u} \\
0 & -{\widetilde{x}}^TG & -x_0 \end{array} \right) \, \, \left( \begin{array}{ccc} 0 & 0 & 0 \\ \widetilde{y} & \widetilde{Y} & 0 \\ 
0 & -{\widetilde{y}}^TG & 0 \end{array} \right) $}}\right] \notag \\
& = \mathrm{tr}\mbox{{\footnotesize $\left( \begin{array}{lcc}
-{\widetilde{u}}^TG\widetilde{y} & * & * \\ * & \widetilde{X}\widetilde{Y} - 
\widetilde{u} \otimes {\widetilde{y}}^TG & * \\ * & * & 0 \end{array} 
\right) $}} \notag \\
& = -2 {\widetilde{u}}^TG\widetilde{y} 
+ \mathrm{tr}\, \widetilde{X}\widetilde{Y}, 
\label{eq-verynewstar} 
\end{align}
for every $\widetilde{y} \in \widetilde{V}$ and every $\widetilde{Y} \in 
\oo (\widetilde{V},G)$. Set $\widetilde{y} =0$ and $\widetilde{Y} = 
{\widetilde{X}}^T$. Then equation (\ref{eq-verynewstar}) reads $0 = 
\mathrm{tr} (\widetilde{X}{\widetilde{X}}^T)$, which implies 
$\widetilde{X} =0$. Now equation (\ref{eq-verynewstar}) reads $0 = {\widetilde{u}}^TG\widetilde{y}$ for every $\widetilde{y} \in 
\widetilde{V}$. But $G$ is invertible. So $\widetilde{u} =0$. Hence 
$X=${\tiny $\begin{pmatrix} x_0 & 0 & 0 \\ \widetilde{x} & 0 & 0 \\
0 & -{\widetilde{x}}^TG & -x_0 \end{pmatrix} $}$=L_{x,e_{n+1}}$, where 
$x = x_0e_0 +\widetilde{x} \in V$. Therefore ${\oo (V,K)}^0_{e_{n+1}} 
\subseteq \{ {\ell }_{L_{v,e_{n+1}}} \in {\oo (V,K)}^{\ast } \setrule \, 
v = v_0e_0 + \widetilde{v} \in V \} $. But 
\begin{displaymath}
\dim {\oo (V,K)}^0_{e_{n+1}} = \dim \oo (V,K) - \dim {\oo (V,K)}_{e_{n+1}} 
=n+1, 
\end{displaymath}
which equals the dimension of the subspace of ${\oo (V,K)}^{\ast }$ 
spanned by covectors of the form ${\ell }_{L_{v,e_{e_{n+1}}}}$ with 
$v = v_0e_0 +\widetilde{v} \in V$. Consequently, equation (\ref{eq-newstar}) 
holds. \hfill $\square $ \medskip 

Now we are in position to prove proposition 9. Supppose 
that the tuples $(V,Y,e_{n+1};K)$ and $(V,Y',e_{n+1};K)$ are 
equivalent. Then there is a $P\in {\Oo (V,K)}_{e_{n+1}}$ and a vector 
$w\in V$ such that $Y'=P(Y+L_{w,e_{n+1}})P^{-1}$. For every 
$Z \in {\oo (V,K)}_{e_{n+1}}$ we have 
\begin{align}
{\ell }_{Y'}(Z) & = {\ell }_{P(Y+L_{w,e_{n+1}})P^{-1}}(Z) = 
{\ell }_{PYP^{-1}}(Z) +{\ell }_{PL_{w,e_{n+1}}P^{-1}}(Z) \notag \\
& = {\ell }_{PYP^{-1}}(Z) +{\ell }_{L_{Pw,e_{n+1}}}(Z), \quad 
     \mbox{since $P\in {\Oo (V,K)}_{e_{n+1}}$} \notag \\
& = {\ell }_{PYP^{-1}}(Z), \quad 
     \mbox{since $Z \in {\oo (V,K)}_{e_{n+1}}$} \notag \\
& = ({\Ad }^T_{P^{-1}}{\ell }_Y)(Z) = {\Ad }^T_{P^{-1}}({\ell }_Y|
    {\oo (V,K)}_{e_{n+1}})(Z). \notag 
\end{align}
So the affine cotype represented by the tuple $(V,Y,e_{n+1};K)$ 
corresponds to the coadjoint orbit of ${\Oo (V,K)}_{e_{n+1}}$ through 
${\ell }_Y|{\oo (V,K)}_{e_{n+1}}$ in 
${\oo (V,K)}^{\ast }_{e_{n+1}}$. Thus the map induced by 
(\ref{eq-newnine}) exists. \medskip 

Suppose that for some $Y,Y' \in \oo (V,K)$ and some 
$P\in {\Oo (V,K)}_{e_{n+1}}$ we have ${\ell }_{Y'}-{\Ad}^T_{P^{-1}}{\ell }_Y 
=0$ on ${\oo (V,K) }_{e_{n+1}}$. In other words, we suppose that 
${\ell}_{Y'}|{\oo (V,K)}_{e_{n+1}}$ lies in the ${\Oo (V,K)}_{e_{n+1}}$ 
coadjoint orbit through ${\ell }_Y|{\oo (V,K)}_{e_{n+1}}$. Then 
${\ell }_{Y'-PYP^{-1}} \in {\oo (V,K)}^0_{e_{n+1}}$. Therefore for 
some $v\in V$ we have ${\ell }_{Y'-PYP^{-1}} = {\ell }_{L_{v,e_{n+1}}}$. 
So 
\begin{displaymath}
Y' = PYP^{-1}+L_{v,e_{n+1}} = P(Y+L_{P^{-1}v, e_{n+1}})P^{-1}.
\end{displaymath}
Hence the tuples $(V,Y,e_{n+1};K)$ and $(V,Y',e_{n+1};K)$ are 
equivalent. Thus the coadjoint orbit of ${\Oo (V,K)}_{e_{n+1}}$ on 
${\oo (V,K)}^{\ast}_{e_{n+1}}$ determines a unique affine cotype. 
Therefore the map induced by (\ref{eq-newnine}) is injective. \medskip 

Since every element of ${\oo (V,K)}^{\ast}_{e_{n+1}}$ is of the form 
${\ell }_Y|{\oo (V,K)}_{e_{n+1}}$ for some $Y\in \oo (V,K)$, the map 
induced by (\ref{eq-newnine}) is surjective. \hfill $\square $ \medskip

Suppose that we are given the affine cotype $\nabla $ with representative
$(\widehat{V},\widehat{Y},$ $\widehat{v};
\widehat{\gamma} )$.
We wish to associate a Gram matrix $K$ to it.
For this, recall that the distinguished type,
represented by $(0,\widehat{V},\widehat{v};\widehat{\gamma })$,
has a representative of the form $(0,V, e_{n+1};K)$, where
$V = \R \times \widetilde{V} \times \R $ and $K=${\tiny $
\begin{pmatrix} 0 & 0 & 1 \\ 0 & G & 0 \\
1 & 0 & 0 \end{pmatrix} $.}
We may replace the representative of the cotype $\nabla $ with
one of the form
$(V,Y, e_{n+1};K)$, where the matrix of $Y$ with respect to the
standard basis $\mathfrak{e}$ is {\tiny $ \begin{pmatrix}
y_0 & -{\widetilde{v}}^{\ast } & 0 \\ \widetilde{y} & \widetilde{Y} &
\widetilde{v} \\ 0 & -{\widetilde{y}}^{\ast } & -y_0 \end{pmatrix} $}. Here $y_0 \in \R $, $\widetilde{v}, \widetilde{y} \in
\widetilde{V}$, $\widetilde{Y} \in \oo (\widetilde{V}, G)$, and 
${\widetilde{v}}^{\ast} = v^TG$. We say that 
the cotype ${\nabla }_{\ell }$, represented by $(\widetilde{V}, 
\widetilde{Y}, \widetilde{v};G)$, is the
\emph{little cotype} of $\nabla $.\footnote{The cotype ${\nabla }_{\ell }$
is called the little cotype because we are imitating the little
subgroup approach of Wigner \cite{wigner} to the representation
theory of the Poincar\'{e} group.} 

\begin{lem} \label{Lem5.2} 
The little cotype ${\nabla }_{\ell }$ does not depend on the choice of representative of the affine cotype $\nabla $. 
\end{lem}

\noindent \textbf{Proof.} Up to isomorphism $(\widetilde{V},G)$ is determined by $\nabla $, so there
is no need to vary $G$ or $K$.
Let $(V, Y, e_{n+1}; K)$ be a representative
of the affine cotype $\nabla $. Suppose that $(V, Y', e_{n+1}; K)$ is another
representative. Then there is a $P \in {\Oo (V,K)}_{e_{n+1}}$ and
a vector $w \in V$ such that
\begin{equation}
Y'= P(Y+L_{w,e_{n+1}})P^{-1} .
\label{eq-5star}
\end{equation}
We now calculate the right hand side of (\ref{eq-5star}) explicitly. With
respect to the standard basis $\mathfrak{e}$ of $(V,K)$, we have
$P=${\tiny $\begin{pmatrix} 1 & 0 & 0 \\
\widetilde{u} & A & 0 \\ -\onehalf {\widetilde{u}}^TG\widetilde{u} &
-{\widetilde{u}}^TGA & 1 \end{pmatrix} $}, where
$\widetilde{u} \in \widetilde{V}$ and $A \in \Oo (\widetilde{V},G)$. 
Therefore
$P^{-1}=${\tiny $ \begin{pmatrix} 1 & 0 & 0 \\
-A^{-1}\widetilde{u} & A^{-1} & 0 \\ -\onehalf {\widetilde{u}}^TG\widetilde{u}
& {\widetilde{u}}^TG & 1 \end{pmatrix} $ }.
Moreover, $Y=$ {\tiny $\begin{pmatrix} y_0 & -{\widetilde{v}}^TG & 0
\\ \widetilde{y} & \widetilde{Y} & \widetilde{v} \\
0 & - {\widetilde{y}}^TG & -y_0 \end{pmatrix} $},
where
$y_0 \in \R $, $\widetilde{v}, \widetilde{y} \in \widetilde{V}$ and
$\widetilde{Y}\in \oo (\widetilde{V}, G)$, and
$L_{w, e_{n+1}} =${\tiny $\begin{pmatrix} w_0 & 0 & 0 \\
\widetilde{w} & 0 & 0 \\ 0 & -{\widetilde{w}}^TG & -w_0 \end{pmatrix} $}, where $w = w_0 e_0 + \widetilde{w}+w_{n+1}e_{n+1}$. So
\begin{eqnarray*}
Y' & = &  P(Y+L_{w,e_{n+1}})P^{-1} \\
& = & \mbox{{\tiny $\left( \begin{array}{ccc} 1 & 0 & 0 \\
\widetilde{u} & A & 0 \\ -\onehalf {\widetilde{u}}^TG\widetilde{u} &
-{\widetilde{u}}^TGA & 1 \end{array} \right) $}}
\mbox{{\tiny $ \left( \begin{array}{ccc} * & -{\widetilde{v}}^TG & 0
\\ {*} & \widetilde{Y} & \widetilde{v} \\
0 & * & * \end{array}  \right) $}}
\mbox{{\tiny $ \left( \begin{array}{ccc} 1 & 0 & 0 \\
-A^{-1}\widetilde{u} & A^{-1} & 0 \\ 
-\onehalf {\widetilde{u}}^TG\widetilde{u}
& {\widetilde{u}}^TG & 1 \end{array} \right) $}} \\
& = & \mbox{{ \tiny $ \left( \begin{array}{ccc} * & -{\widetilde{v}}^TG
& 0 \\
{*} &
-\widetilde{u} \otimes {\widetilde{v}}^TG + A\widetilde{Y} & 
A \widetilde{v} \\
{*} &
{ * }
& *
\end{array} \right) $ }}
\mbox{{\tiny $ \left( \begin{array}{ccc} 1 & 0 & 0 \\
-A^{-1}\widetilde{u} & A^{-1} & 0 \\ -\onehalf 
{\widetilde{u}}^TG\widetilde{u}
& {\widetilde{u}}^TG & 1 \end{array} \right) $}} \\
& = & \mbox{{\tiny $ \left( \begin{array}{ccc} b_0 &
 -{ (\widetilde{v})'}^TG &
0 \\
\widetilde{b} & {\widetilde{Y}}' & {\widetilde{v}}' \\
0 & -{\widetilde{b}}^T G & -b_0 \end{array} \right) ,$ }}
\end{eqnarray*}
where $b_0 \in \R $, $\widetilde{b} \in \widetilde{V}$,
\begin{eqnarray*}
{\widetilde{Y}}' & = & A\widetilde{Y}A^{-1} - \widetilde{u} \otimes
{\widetilde{v}}^TGA^{-1} + A\widetilde{v} \otimes {\widetilde{u}}^TG \\
& = & A\widetilde{Y}A^{-1} - \widetilde{u} \otimes 
(A\widetilde{v})^{\ast } +
A\widetilde{v} \otimes {\widetilde{u}}^{\ast } \\
& = &  A\widetilde{Y}A^{-1} + L_{-\widetilde{u}, A\widetilde{v}},
\end{eqnarray*}
and ${\widetilde{v}}' = A\widetilde{v}$. Thus the little cotype
${\nabla}_{\ell }$, as computed from $(V, Y', e_{n+1}; K)$,
is represented by the tuple
$(\widetilde{V}, {\widetilde{Y}}', {\widetilde{v}}'; G)$, which does
not depend on the vector $w$. Since ${\widetilde{Y}}' = A(\widetilde{Y} +
L_{-A^{-1}\widetilde{u}, \widetilde{v}})A^{-1}$ 
and ${\widetilde{v}}' = A\widetilde{v}$, the tuple
$(\widetilde{V}, {\widetilde{Y}}', {\widetilde{v}}'; G)$ is equivalent
to the tuple $(\widetilde{V}, \widetilde{Y}, \widetilde{v}; G)$. But this tuple
depends only on the representative $(V, Y, e_{n+1};K)$ and
\emph{not} the representative $(V, Y', e_{n+1};K)$ of
the cotype $\nabla $. So the little cotype ${\nabla }_{\ell }$ does
not depend on the choice of representative of the affine cotype $\nabla $.
\hfill $\Box $ 

\begin{lem} \label{Lemma 5.3}
Let $\nabla $ be an affine
cotype. Then $\nabla $ is uniquely determined by its little
cotype ${\nabla }_{\ell }$.
\end{lem} 

\noindent \textbf{Proof.} Suppose that the affine cotypes
$\nabla $ and ${\nabla }'$, represented by the tuples
$(V,Y,e_{n+1};K) $ and $(V,Y',e_{n+1}; K) $, both have 
the little cotype ${\nabla }_{\ell }$.
Say $Y=${\tiny $ \begin{pmatrix}
u_0 & -{\widetilde{w}}^{\ast} & 0 \\ 
\widetilde{u} & \check{Y} & \widetilde{w} \\
0 & -{\widetilde{u}}^{\ast} & -u_0 \end{pmatrix} $} and
$Y'=${\tiny $\begin{pmatrix}
u'_0 & -({\widetilde{w}}')^{\ast} & 0 \\ 
\widetilde{u}' & (\check{Y})' & {\widetilde{w}}' \\
0 & -({\widetilde{u}}')^{\ast} & 0 \end{pmatrix} $}, 
where $u_0,u'_0 \in \R $, $\widetilde{u},{\widetilde{u}}', \widetilde{w}, 
{\widetilde{u}}' \in \widetilde{V}$, and
$\check{Y}, (\check{Y})' \in \oo (\widetilde{V},G)$.
Thus ${\nabla }_{\ell }$ is represented by the tuples 
$(\widetilde{V}, \check{Y}, \widetilde{w};G)$ and 
$(\widetilde{V}, \check{Y}', \widetilde{w}';G)$, which are equivalent. In other words, there is a $\widetilde{A} \in \Oo (\widetilde{V},G)$ and a vector
$\widetilde{u} \in \widetilde{V}$ such that $\widetilde{A}\widetilde{w} =
{\widetilde{w}}'$ and
\begin{displaymath}
({\check{Y}})' = \widetilde{A}(\check{Y} +L_{\widetilde{w},\widetilde{u}}){\widetilde{A}}^{-1}.
\end{displaymath}
Let $A=${\tiny $ \begin{pmatrix} 1 & 0 & 0 \\
\widetilde{A}\widetilde{u} & \widetilde{A} & 0 \\
-\onehalf {\widetilde{u}}^TG\widetilde{u} & 
-{\widetilde{u}}^{\ast} & 1 \end{pmatrix} $}. Then 
$A\in {\Oo (V,K)}_{e_{n+1}}$. Now
\begin{align}
AYA^{-1} & =   
\mbox{{\footnotesize $\left( \begin{array}{crc} 
1                            &      0       & 0 \\
\widetilde{A} \widetilde{u}  &\widetilde{A} & 0 \\
-\onehalf {\widetilde{u}}^TG\widetilde{u} & -{\widetilde{u}}^{\ast } & 1 
\end{array} \right) $}} 
\, 
\mbox{{\footnotesize $\left( \begin{array}{crr} 
u_0 & -{\widetilde{w}}^{\ast} & 0 \\
\widetilde{u} & \check{Y} & \widetilde{w} \\
0 & -{\widetilde{u}}^{\ast} & -u_0 \end{array} \right)$}} \, 
\mbox{{\footnotesize $\left( \begin{array}{ccc} 
1 & 0 & 0 \\
-\widetilde{u} & {\widetilde{A}}^{-1} & 0 \\
-\onehalf {\widetilde{u}}^TG\widetilde{u} & {\widetilde{u}}^{\ast} 
{\widetilde{A}}^{-1} & 1 \end{array} \right) $}} \notag \\ 
& = \mbox{{\footnotesize $\left( \begin{array}{crr} 
v_0 & - {\widetilde{r}}^{\ast } & 0 \\
\widetilde{v} & \check{Z} & \widetilde{r} \\
0 & -{\widetilde{v}}^{\ast } & -v_0 \end{array} \right) $},} \notag
\end{align}
where 
\begin{align}
v_0 & = u_0 +{\widetilde{w}}^{\ast}(\widetilde{u}) \notag \\
\widetilde{v} & = u_0\widetilde{A}\widetilde{u} + \widetilde{A}\widetilde{u} 
\otimes {\widetilde{w}}^{\ast} (\widetilde{u}) 
-\widetilde{A}\check{Y}\widetilde{u} 
-\onehalf ({\widetilde{u}}^TG\widetilde{u})\widetilde{A}\widetilde{w} 
\notag \\
\widetilde{r} & = \widetilde{A}\widetilde{w} = {\widetilde{w}}' \notag \\
\check{Z} & = -\widetilde{A}\widetilde{u} \otimes {\widetilde{w}}^{\ast } 
{\widetilde{A}}^{-1} + \widetilde{A} \check{Y}{\widetilde{A}}^{-1} +
\widetilde{A}\widetilde{w} \otimes {\widetilde{u}}^{\ast} {\widetilde{A}}^{-1} 
\notag \\
& = -\widetilde{A}\widetilde{u} \otimes (\widetilde{A}\widetilde{w})^{\ast } 
+\widetilde{A}\widetilde{w} \otimes (\widetilde{A}\widetilde{u})^{\ast } +
\widetilde{A}\check{Y}{\widetilde{A}}^{-1} \notag \\
& = \widetilde{A}\check{Y}{\widetilde{A}}^{-1} +L_{\widetilde{A}\widetilde{w}, 
\widetilde{A}\widetilde{u}} = {\widetilde{A}}(\check{Y} + L_{\widetilde{w}, 
\widetilde{u}}){\widetilde{A}}^{-1} = (\check{Y})'. \notag 
\end{align}

So 
\begin{align}
AYA^{-1} & = \mbox{{\footnotesize $\left( \begin{array}{ccc} 
0 & -({\widetilde{w}}')^{\ast} & 0 \\
0 & (\check{Y})' & {\widetilde{w}}' \\
0 & 0 & 0 \end{array} \right) $}} + L_{v,e_{n+1}}, \quad 
\mbox{where $v = v_0e_0 + \widetilde{v} \in V$} \notag  \\
& = \mbox{{\footnotesize $\left( \begin{array}{crc}
u'_0 & -({\widetilde{w}}')^{\ast } & 0 \\
{\widetilde{u}}' & (\check{Y})' & {\widetilde{w}}' \\
0 & -({\widetilde{u}}')^{\ast} & -u'_0 \end{array} 
\right) + L_{v-u',e_{n+1}} $}}, \quad 
\mbox{where $u' = u'_0 e_0 + {\widetilde{u}}' \in V$} \notag \\
& = Y' +L_{v-u',e_{n+1}}, \notag 
\end{align}
which implies 
\begin{displaymath}
Y' = A(Y+L_{A^{-1}(u'-v), e_{n+1}})A^{-1}.
\end{displaymath}
In other words, the tuples $(V,Y,e_{n+1};K)$ and $(V,Y',e_{n+1};K)$ are 
equivalent. Thus the affine cotypes $\nabla $ and ${\nabla}'$ are 
equal. \hfill $\square $   

\begin{remark} 
Given a  cotype ${\nabla }_{\ell }$
it is very easy to construct a cotype ${\nabla }$
having ${\nabla }_{\ell }$ as little cotype.
Indeed if $(\widetilde{V}, \widetilde{Y}, \widetilde{v};G)$
represents ${\nabla }_{\ell }$, one forms $V$, $K$ in the usual way and
takes a representative of the form
$(V,Y, e_{n+1};K)$, where the matrix of $Y$ with respect to the
standard basis $\mathfrak{e}$ is {\tiny $\begin{pmatrix}
0 & -{\widetilde{v}}^{\ast } & 0 \\ 0 & \widetilde{Y} &
\widetilde{v} \\ 0 & 0 & 0 \end{pmatrix} $}. 
\end{remark}

The following proposition follows immediately from the above. 

\begin{prop} \label{prop5.4} 
There is a bijection between little cotypes and coadjoint orbits. 
\end{prop} 

Let $\nabla $ be a cotype represented by
the tuple $(V,Y,v; \gamma )$. Suppose
that $V = V_1 \oplus V_2$, where $V_i$ are $Y$-invariant,
$\gamma $-nondegenerate and $\gamma $-orthogonal subspaces
such that $V_2\neq\{ 0 \} $ and $v \in V_1$. Then we say that 
$\nabla $ is a
\emph{sum} of the cotype $\widetilde{\nabla}$, represented by the 
tuple $(V_1,Y|V_1, v;\gamma |V_1)$, and a \emph{type} $\Delta $, represented by
$( Y|V_2,V_2; {\gamma }|V_2)$. We write $\nabla = \widetilde{\nabla } +
\Delta $. If $V_1 = \{ 0 \} $, then $v =0$ and $\widetilde{\nabla }$ is
the \emph{zero cotype}, represented by the tuple $(\{ 0 \}, 0,0;0)$ and
denoted by $\mathbf{0}$. We say the a cotype is \emph{indecomposable} if it
cannot be written as the sum of a cotype and a type. A nonzero cotype
$\nabla $, represented by the tuple $(V, Y, v; \gamma )$ is decomposable
if there is a proper, $Y$-invariant subspace of $V$, which
contains the vector $v$ and on which $\gamma $ is nondegenerate.
Conversely, if $\nabla $ is decomposable, then there is a representative
$(V, Y, v; \gamma )$ so that there is
a proper, $Y$-invariant subspace of $V$, which
contains the vector $v$ and on which $\gamma $ is nondegenerate.
Let us call such a representative \emph{adapted} to the decomposition.

\begin{lem}\label{Lemma 5.6} 
Every cotype, 
which is \emph{not} affine, is
the sum of a unique indecomposable cotype, which is either the zero
cotype or a nonzero $1$-dimensional cotype, and a type.
\end{lem} 

\noindent \textbf{Proof.} Let $(V,Y,v; \gamma )$ represent the nonaffine
cotype $\nabla $. Suppose that $v=0$. Write $V = \{ 0 \} \oplus V$.
Then $\{ 0 \} $ and $V$ are $Y$-invariant, $\gamma $-orthogonal,
and $\gamma $-nondegenerate. Hence $\nabla $ is the sum of the zero
cotype $\mathbf{0}$ and a type $\Delta $, represented by $(Y,V; \gamma )$.
Now suppose that $v \ne 0$. Because $\nabla $ is not affine, $v$ is
not $\gamma $-isotropic, that is, $\gamma (v,v) = \varepsilon \, {\alpha }^2$,
where ${\varepsilon }^2 =1$ and $\alpha >0$. Since $\spann \{ v \} $ is
$\gamma $-nondegenerate, its orthogonal complement $\widetilde{V} =
{\spann \{ v \} }^{\gamma }$ is also $\gamma $-nondegenerate. Let
$\widetilde{\mathfrak{f}} = \{ e_1, \ldots , e_n \}$ be a
basis of $\widetilde{V}$ such that the matrix of $\widetilde{\gamma } =
\gamma |\widetilde{V}$ is $F$. Then $\mathfrak{f} = \{ e_1, \ldots , e_n,
e_{n+1} = v\} $ is a basis of $V$ such that
the matrix of $\gamma $ with respect to $\mathfrak{f}$ is $G=${\tiny
$\begin{pmatrix} F & 0 \\ 0 & \varepsilon {\alpha }^2
\end{pmatrix} $.} Since $Y \in \oo (V, \gamma )$, the matrix of
$Y$ with respect to the basis $\mathfrak{f}$ is $Y=${\tiny $
\begin{pmatrix} \widetilde{Y} & \varepsilon {\alpha }^2 \, \widetilde{v} \\
-{\widetilde{v}}^TF & 0 \end{pmatrix} $,} where $\widetilde{Y} \in
\oo (\widetilde{V}, \widetilde{\gamma })$ and $\widetilde{v} \in
\widetilde{V}$.
Thus the tuple $(V,Y, e_{n+1}; G)$ represents the cotype $\nabla $. For
every $w = \widetilde{w} + w_{n+1}e_{n+1} \in \widetilde{V} \oplus
\spann \{ e_{n+1} \} $, the matrix of $L_{w, e_{n+1}}$ with respect
to the basis $\mathfrak{f}$ is {\tiny $\begin{pmatrix}
0 & \varepsilon {\alpha }^2 \, \widetilde{w} \\ -{\widetilde{w}}^TF & 0
\end{pmatrix} $}$\in \oo (V,G)$, since
\begin{eqnarray*}
L_{w,e_{n+1}}(e_i) & = & e^{\ast }_{n+1}(e_i)w -w^{\ast }(e_i)e_{n+1} \\
& = & - (w^TGe_i)e_{n+1} \, = - ({\widetilde{w}}^TFe_i)e_{n+1}, \qquad
\mbox{for $1 \le i \le n$} \\
L_{w,e_{n+1}}(e_{n+1}) & = & e^{\ast }_{n+1}(e_{n+1})w
-w^{\ast }(e_{n+1})e_{n+1} \\
& = & (e^T_{n+1}Ge_{n+1})w - (w^TGe_{n+1})e_{n+1} \\
& = & \varepsilon {\alpha }^2\,  (w-w_{n+1}e_{n+1}) \, = \,
\varepsilon {\alpha }^2 \, \widetilde{w}.
\end{eqnarray*}
Therefore we may write $Y = ${\tiny $\begin{pmatrix} 
\widetilde{Y} & 0 \\ 0 & 0
\end{pmatrix} $}$+ L_{\varepsilon {\alpha }^2 \, \widetilde{v},
e_{n+1}}$, which implies that the tuple $(V,Y, e_{n+1}; G)$ is 
equivalent to the
tuple $(V, \check{Y} =$ 
{\tiny $\begin{pmatrix} \widetilde{Y} & 0 \\ 0 & 0
\end{pmatrix} $}, $e_{n+1};G)$. Now the
subspace $\spann \{ e_{n+1} \} $ is $G$-nondegenerate, since
the matrix of $G$ restricted to $\spann \{ e_{n+1} \} $ is
$(\varepsilon {\alpha }^2)$, which is nonzero. From
$\check{Y}e_{n+1} =0$, it follows that
$\spann \{ e_{n+1} \} $ is $\check{Y}$-invariant. Clearly, the space
$\widetilde{V} = {\spann \{ e_{n+1} \} }^G$ is also $\check{Y}$-invariant.
Therefore
the cotype $\nabla $, represented by the tuple $(V, \check{Y}, e_{n+1};G)$,
is the sum of a $1$-dimensional cotype $\widetilde{\nabla}$, represented by
the tuple $(\spann \{ e_{n+1} \},$ $ 0, e_{n+1}; 
(\varepsilon {\alpha }^2))$, and
a type $\Delta $, represented by $( \widetilde{Y},\widetilde{V}; F)$.
\hfill $\Box $ 

\begin{lem} \label{Lemma 5.7} 
Every affine cotype can be written as a sum of an indecomposable affine cotype 
and a sum of indecomposable types. This decomposition is unique up to reordering of the summands which are types. 
\end{lem}

\noindent \textbf{Proof.} Suppose that we are given an affine cotype
$\nabla $. Then $\nabla $ is uniquely determined by its little
cotype ${\nabla }_{\ell }$, where $\dim {\nabla }_{\ell } < \dim \nabla $.
This correspondence respects decomposition:
if ${\nabla }_{\ell }$ is decomposable, then reconstructing
${\nabla }$ as in the remark above, one finds that ${\nabla }$
is decomposable. Conversely, if ${\nabla }$
is decomposable, then using a representative adapted to a decomposition
one finds that ${\nabla }_{\ell }$ is decomposable.
If ${\nabla }_{\ell }$ is again affine, we look at its little cotype.
Repeating this process a finite number of times, we obtain either the
zero cotype and we stop or we obtain a nonzero cotype $\widehat{\nabla }$
which is not affine. By lemma \ref{Lemma 5.6}
$\widehat{\nabla }$ is a unique sum
of a cotype $\widetilde{\nabla }$, which is either the zero cotype or
a nonzero $1$-dimensional cotype and a type $\Delta $. By results of
\cite{burgoyne-cushman}, the type $\Delta $ is a sum of indecomposable
types, which is unique up to reordering the summands. This completes
the proof of the lemma. \hfill $\Box $ \medskip

We now classify indecomposable affine tuples. Let $(V, Y, e_{n+1}; K )$ be an indecomposable affine tuple with respect to the standard orthogonal basis $\mathfrak{e} = 
\{ e_0, e_1, \ldots , e_n,$ $e_{n+1} \} $ of $V$. The matrix of 
$Y$ with respect to $\mathfrak{e} $ is {\tiny $\begin{pmatrix} 
x_0 & -{\widetilde{y}}^{\, T}G & 0 \\
\widetilde{x} & \widetilde{Y} & \widetilde{y} \\
0 & -{\widetilde{x}}^{\, T}G & - x_0  \end{pmatrix} $}, where $x_0 \in \R$; 
$\widetilde{x}, \widetilde{y} \in \widetilde{V} = \spann \{ e_1, \ldots , e_n \} $; 
and $\widetilde{Y}\in \oo (\widetilde{V}, G)$. 
Let $x = x_0e_0 + \widetilde{x} \in V$. The matrix 
of $L_{x,e_{n+1}}$ with \linebreak 
respect to $\mathfrak{e}$ is 
{\tiny $\begin{pmatrix} x_0 & 0 & 0 \\ \widetilde{x} & 0 & 0 \\
0 & -{\widetilde{x}}^{\, T}G & -x_0 \end{pmatrix}$}. 
Consider the tuple $(V, Y', e_0; K )$.  
Here $Y' = PY''P^{-1}$ with $P \in \Oo (V, K)$ given by $Pe_0 = e_{n+1}$, $P_{|\widetilde{V}} = 
{\mathrm{id}}_{\widetilde{V}}$ and $Pe_{n+1} = e_0$. Here $Y'' = Y + L_{-x,e_{n+1}} = 
${\tiny $\begin{pmatrix}
0 & -{\widetilde{y}}^{\, T}G & 0 \\ 0 & \widetilde{Y} & \widetilde{y} \\ 
0 & 0 & 0 \end{pmatrix}$}. Since the matrix of $P$ with respect to 
$\mathfrak{e}$ is {\tiny $\begin{pmatrix} 0 & 0 & 1 \\ 0 & 
{\mathrm{id}}_{\widetilde{V}} & 0 \\ 
1 & 0 & 0 \end{pmatrix}$}, we get $Y' =${\tiny $\begin{pmatrix} 0 & 0 & 0 \\\widetilde{y} & \widetilde{Y} & 0 \\ 0 & -{\widetilde{y}}^{\, T}G & 0 
\end{pmatrix}$}. Thus $Y'$ is the matrix of $Y''$ with respect to the orthogonal basis $P\mathfrak{e}  = \{ e_{n+1}, e_1, \ldots , e_n, e_0 \}$. The 
tuple $(V, Y', e_0; K )$ is affine, is equivalent to the tuple 
$(V, Y, e_{n+1}; K)$, and has $Y'e_0 =0$. Suppose that 
the tuple $(V, Y', e_0; K)$ is indecomposable. Let $\widehat{V}$ be the generalized eigenspace of $Y'$ corresponding to the eigenvalue $0$. 
Then $e_0 \in \widehat{V}$ and $\widehat{V}$ is a $Y'$ invariant subspace on which $Y'$ is nilpotent and $K|\widehat{V}$ is nondegenerate. Since 
$(V,Y', e_0; K)$ is indecomposable, it follows that $\widehat{V} = V$. Thus 
$(V, Y', e_0; K)$ is an indecomposable, nilpotent, affine tuple. 
Hence $(V, Y', e_0; K)$ is a triple, which represents an indecomposable, nilpotent distinguished type. Such distinguished types were classified in 
proposition 5, which we restate as   

\begin{prop}\label{Proposition 5.8} 
Let $\nabla$ be an indecomposable affine cotype of dimension $n$, 
which is represented by the nilpotent affine tuple 
$(V,Y, v^0; K)$. Then exactly one of the following alternatives holds.
\begin{itemize}
\item[$\mathbf{1}$.] $\dim V = n$ is \emph{even}, say $n=2h+2$, $h\geq0$.
There is a representative $(V,Y,v^0;K)$ of $\nabla$ such that the
following holds. There is an orthogonal basis ${\epsilon }_{2(h+1)}$ 
of $V$ given by  
\begin{displaymath}
\{ (-1)^hz, (-1)^{h-1}Yz, \, \ldots ,\,  -Y^{h-1}z,  Y^hz\,   ; 
w, \, Yw, \ldots , Y^hw = v^0 \}
\end{displaymath}
with $Y^{h+1} =0$ such that $K (Y^jz, Y^{h-j}w) = 
(-1)^{h-j}$ for $j =0, 1, \ldots , h$ and $K (Y^kz, Y^{\ell }w) = 0$ if 
$k + \ell \ne h$. In other words, the matrix of $Y$ with respect to the 
basis ${\mathfrak{e}}_{2(h+1)}$ is 
\begin{displaymath} 
\mbox{\footnotesize ${\mathcal{N}}_{2(h+1)} = 
\begin{pmatrix} -N_{h+1} & 0 \\ 0 & N_{h+1} \end{pmatrix}$} 
= \mbox{\footnotesize $ \left( \begin{array}{r|cc|c} 
0 & 0 & 0 & 0\\ \hline 
-e_1 & -N_h & 0 & 0 \\
0  &  0 & N_h & 0 \\ \hline
\rule{0pt}{10pt} 0 & 0 & e^T_h & 0 \end{array} \right) $,} 
\end{displaymath}%
where $N_{h+1}$ is an $(h+1) \times (h+1)$ lower Jordan block with 
$N_1 = 0$. The matrix of $K$ with respect ${\mathfrak{e}}_{2(h+1)}$ is $K_{2(h+1)} =${\tiny $\left( \begin{array}{c|c|c}0 & & 1 \\ \hline 
 & K_{2h} & \\ \hline
1 &  & 0 \end{array} \right) $,} 
where $K_2=${\tiny $\begin{pmatrix} 0 & 1 \\ 1 & 0 \end{pmatrix}$} 
and the index of $K_{2(h+1)}$ is $h+1$. There is no modulus. We use the notation ${\nabla}_n(0,0)$ for the cotype $\nabla $.
\item[$\mathbf{2}$.] $\dim V = n$ is \emph{odd}, say $n = 2h+1$, 
$h \ge 0$. There is a representative $(V,Y,v^0; K)$ of $\nabla$ such that the following holds. There is an orthogonal basis ${\mathfrak{e}}_{2h+1}$ of $V$ given by
\begin{align*}
&  \{ w, \, Yw, \, \ldots , Y^{h-1}w; 
\varepsilon Y^hw; \notag \\
& \hspace{.5in} \, (-1)^{h+1}{\varepsilon }Y^{h+1}w, 
\, (-1)^{h+2}\varepsilon Y^{h+2}w, \, \ldots , \varepsilon Y^{2h}w = v^0 \}
\end{align*}
with $Y^{2h+1} =0$ such that $K(Y^jw, Y^{2h-j}w) = (-1)^j\varepsilon $ for $j=0,1, \ldots , h$ and $K(Y^kw, Y^{\ell }w) $ $= 0 $ if $k+\ell \ne 2h$. Here 
${\varepsilon }^2 =1$. In other words, the matrix of $Y$ with respect to 
the ${\mathfrak{e}}_{2h+1}$ basis is 
\begin{displaymath}
\mbox{\footnotesize $ {\mathcal{N}}_{2h+1} = \left( \begin{array}{c|c|c} 
N_h & & \\ \hline 
\varepsilon & 0 & \\ \hline 
& & -N_h \end{array} \right) $,}
\end{displaymath}
where $N_1 =0$; while the matrix of $K$ is $K_{2h+1}=${\tiny 
$ \left( \begin{array}{c|c|c} 0 & & 1 \\ \hline 
& K_{2h-1} & \\ \hline
1 & & 0 \end{array} \right) $}, where 
$K_1 = \big( (-1)^h\varepsilon \big) $. The index of $K_{2h+1}$ 
is {\tiny $\left\{ \begin{array}{ll} h+1, & \mbox{if $\varepsilon = 1$} \\
h+2, & \mbox{if $\varepsilon = -1$.} \end{array} \right. $} There is a modulus 
$\mu > 0$, where $v^0 = \mu Y^{2h}w$. We use the notation 
${\nabla }^{\varepsilon }_n(0), \, \mu $ for the cotype $\nabla $.
\end{itemize}
\end{prop} 

\noindent \textbf{Proof.} The existence of the bases in cases 1 and 2 follow from 
the proof of proposition 5. In case 1 the tuple representing the little cotype 
${\nabla }_{\ell }$ corresponding to the nilpotent affine cotype 
$\nabla = {\nabla }_n(0,0)$, represented by the tuple 
$({\R }^{2(h+1)}, {\mathcal{N}}_{2(h+1)}, e_{2(h+1)};$ $K_{2(h+1)})$, is 
the nilpotent affine tuple 
$({\R }^{2h}, {\mathcal{N}}_{2h}, e_{2h};$ $K_{2h})$, which represents 
the indecomposable cotype ${\nabla }_{n-2}(0,0)$. 
In case 2 the tuple representing the little cotype ${\nabla }_{\ell }$
corresponding to the nilpotent affine cotype $\nabla = 
{\nabla }^{\varepsilon}_n(0), \, \mu $, represented by the tuple 
$({\R }^{2h+1}, {\mathcal{N}}_{2h+1}, $ $\mu e_{2h+1}; K_{2h+1})$, is 
the nilpotent affine tuple 
$({\R }^{2h-1}, {\mathcal{N}}_{2h-1}, \mu e_{2h-1}; K_{2h-1})$, which 
represents the cotype ${\nabla }^{\varepsilon }_{n-2}(0), \, \mu $.  
\hfill $\square $ 

\begin{remark} It is noteworthy that we can choose the representatives in proposition \ref{Proposition 5.8} to have nilpotent $Y$.
\end{remark} 

\begin{remark} (The curious bijection) 
There is a curious bijection between the representatives that we choose here for indecomposable affine cotypes and the representatives that we used for indecomposable distinguished types in proposition \ref{Proposition 2.2}.
The bijection preserves dimension, index, modulus, and Jordan type. It follows that we also get a bijection between affine cotypes and distinguished types with the same underlying $(V;\gamma)$.
In other words, we get a bijection between adjoint orbits and coadjoint orbits for any affine orthogonal group. 
\end{remark}

\section{Coadjoint orbits of the Poincar\'{e} group}\label{section 6}

In this section we use the theory of \S \ref{section 5} to classify the
coadjoint orbits of the Poincar\'{e} group ${\Oo (4,2)}_{e_5}$. \bigskip

\noindent \begin{tabular}{lccc|lccc}
&\multicolumn{1}{l}{affine cotype} & \multicolumn{1}{c}{dim} &
\multicolumn{1}{c}{index} & & \multicolumn{1}{l}{affine cotype} &
\multicolumn{1}{c}{dim} & \multicolumn{1}{c}{index} \\ \hline
$1.$ & ${\nabla }^{-}_5(0), \, \mu $ & $5$ & $3$ &
$4.$ & ${\nabla }^{+}_3(0), \, \mu $ & $3$ & $1$ \\
$2.$ & ${\nabla }_4(0,0)$ & $4$ & $2$ &
$5.$ & ${\nabla }_2(0,0)$ & $2$ & $1$ \\
$3.$ & ${\nabla }^{-}_3(0), \, \mu $ & $3$ & $2$ &
\end{tabular}\bigskip 

\noindent \hspace{.6in}Table 4. Possible $\oo (V,K)$-indecomposable affine cotypes. \bigskip 

Let $(V, \gamma )$ be a real vector space with a nondegenerate
inner product $\gamma $ of signature $(m,p) = (4,2)$. Suppose
that the tuple $(V,Y',v;\gamma )$ represents an affine cotype in
$\Oo (V, \gamma )$. Since $\Oo (V, \gamma )$ acts transitively on
the collection of nonzero $\gamma $-isotropic vectors in $V$, there
is a $P \in \Oo (V, \gamma )$ such that $Pv =e_5$. Hence the
tuple $(V, Y=PY'P^{-1}, e_5; \gamma )$ is equivalent to
$(V,Y,v;\gamma )$. Because $e_5$ is $\gamma $-isotropic and
$\gamma $ is nondegenerate on $V$, there is a $\gamma $-isotropic
vector $e_0 \in V$ such that $\gamma (e_0,e_5) =1$. In other words,
$H = \spann \{ e_0, e_5 \} $ is a hyperbolic plane in $V$. Because
$\gamma |H$ is nondegenerate, we can extend $\{ e_0, e_5 \} $ to
a $\gamma $-orthonormal basis $\mathfrak{e} = \{ e_0, e_1,
\ldots , e_4, e_5 \} $ of $V$ such that the matrix of $\gamma $ with
respect to $\mathfrak{e}$ is $K=${\tiny $\begin{pmatrix}
0 & 0 & 1 \\ 0 & G & 0 \\ 1 & 0 & 0 \end{pmatrix} $,} where
$G^T = G$, $G^2 = I$, and $G$ has signature $(1,3)$. Thus using the basis
$\mathfrak{e}$ the tuple $(V,Y,e_5; \gamma )$ is the tuple
$(V,Y, e_5; K)$. \bigskip 

\noindent \hspace{.5in}\begin{tabular}{lccc|lccc}
&\multicolumn{1}{c}{type} & \multicolumn{1}{c}{dim} &
\multicolumn{1}{c}{index} & & \multicolumn{1}{c}{type} &
\multicolumn{1}{c}{dim} & \multicolumn{1}{c}{index} \\ \hline
1. &${\Delta }_1(\zeta , \mathrm{RP})$ & $4$ & $2$ &
6. & ${\Delta }^{-}_0(\zeta , \mathrm{IP})$ & $2$ & $2$ \\
2. & ${\Delta }^{\varepsilon }_1(\zeta ,\mathrm{IP})$ & $4$ & $2$ &
7. & ${\Delta }_0(\zeta ,\mathrm{RP})$ & $2$ & $1$ \\
3. & ${\Delta }_1(0,0)$ & $4$ & $2$ &
8. & ${\Delta }^{+}_0(\zeta , \mathrm{IP})$ & $2$ & $0$ \\
4. & ${\Delta }^{+}_2(0)$ & $3$ & $2$ &
9. & ${\Delta }^{-}_0(0)$ & $1$ & $1$ \\
5. & ${\Delta }^{-}_2(0)$ & $3$ & $1$ &
10. & ${\Delta }^{+}_0(0)$ & $1$ & $0$
\end{tabular}  \medskip

\noindent Table 5. Possible $\oo (V,K)$-indecomposable, which 
appear as a summand in the type $\Delta $. \bigskip

Without loss of generality we can begin with an \emph{affine} cotype
$\nabla $ in $\oo (V, K)$ represented by the tuple
\begin{equation}
({\R }^6, Y=\mbox{{\tiny $ \left( \begin{array}{ccc} y_0
& -{\widetilde{x}}^tG & 0 \\
\widetilde{y} & \widetilde{Y} & \widetilde{x} \\ 0 & -{\widetilde{y}}^TG &
-y_0 \end{array} \right) $}}, e_5; K=\mbox{{\tiny $ \left(
\begin{array}{ccc} 0 & 0 & 1 \\ 0 & G & 0 \\ 1 & 0 & 0 \end{array}
\right) $}}),
\label{eq-s6one}
\end{equation}
where $y_0 \in \R $, $\widetilde{x}, \widetilde{y} \in {\R }^4$, and
${\widetilde{Y}}^TG + G\widetilde{Y} =0$,
that is, $\widetilde{Y} \in \oo ({\R }^4, G)$. By proposition
\ref{Lemma 5.7} we can write
$\nabla = \widetilde{\nabla } + \Delta $, where
the possible indecomposable affine cotypes $\widetilde{\nabla }$
in $\oo (V,K)$ are listed in table 4, and the possible
indecomposable summands of the $\oo (V,K)$ type $\Delta $ are listed
in table 5.

Therefore the possible decompositions of the affine cotype $\nabla $
into a sum of an indecomposable affine cotype $\widetilde{\nabla }$ and
a sum of indecomposable types is given in table 6.  

\begin{center}
\begin{tabular}{rlcc}
& \multicolumn{1}{c}{indecomposable affine cotypes} & & \\
& \multicolumn{1}{c}{and sum of indecomposable types} &
\multicolumn{1}{c}{dim} & \multicolumn{1}{c}{index} \\ \hline
$1$. &  ${\nabla }^{-}_5(0), \, \mu +{\Delta }^{-}_0(0) $
& $5+1$ & $3+1$  \\ \hline
$2$. & ${\nabla }_4(0,0) + {\Delta }^{-}_0(\zeta , \mathrm{IP})$
&  $4+2$ & $2+2$  \\
$3$. & ${\nabla }_4(0,0) +{\Delta }^{-}_0(0)+ {\Delta }^{-}_0(0)$
& $4+2$ & $2+2$ \\ \hline
$4$. &  ${\nabla }^{-}_3(0), \, \mu \, +{\Delta }^{+}_2(0)$
&  $3+3$ & $2+2$  \\
$5$. & ${\nabla }^{-}_3(0), \, \mu \, +{\Delta }^{-}_0(\zeta , 
\mathrm{IP}) +
{\Delta }^{+}_0(0)$ & $3+3$ & $2+ 2$ \\
$6$. &  ${\nabla }^{-}_3(0), \, \mu \, + {\Delta }_0(\zeta , \mathrm{RP})
+ {\Delta }^{-}_0(0)$
& $3+3$ & $2+2$ \\
$7$. & ${\nabla }^{-}_3(0), \, \mu \, + {\Delta }^{-}_0(0)+
{\Delta }^{-}_0(0)+
{\Delta }^{+}_0(0)$
& $3+3$ & $2+2$ \\  \hline
$8$. &  ${\nabla }^{+}_3(0), \, \mu \, +{\Delta }^{-}_0(\zeta , 
\mathrm{IP}) +
{\Delta }^{-}_0(0)$
&  $3+3$ & $1+3$  \\
$9$. & ${\nabla }^{+}_3(0), \, \mu \, +  {\Delta }^{-}_0(0)+
{\Delta }^{-}_0(0)+
{\Delta }^{-}_0(0)$
& $3+3$ & $1+ 3$ \\ \hline
$10$. &   ${\nabla }_2(0,0)+{\Delta }^{+}_2(0) + {\Delta }^{-}_0(0)$ &
$2+4$ & $1+3$ \\
$11$. & ${\nabla }_2(0,0) +{\Delta }^{-}_0(\zeta , \mathrm{IP})+
{\Delta }_0(\zeta , \mathrm{RP}) $ & $2+4$ & $1+3$ \\
$12$. &$ {\nabla }_2(0,0)+{\Delta }^{-}_0(\zeta , \mathrm{IP})+
{\Delta }^{-}_0(0)+ {\Delta }^{+}_0(0)$ & $2+4$ & $1+3$ \\
$13$. &${\nabla }_2(0,0)+{\Delta }_0(\zeta , \mathrm{RP})+
{\Delta }^{-}_0(0)+
{\Delta }^{-}_0(0)$ & $2+4$ & $1+3$ \\
$14$. & ${\nabla }_2(0,0) +{\Delta }^{-}_0(0)+{\Delta }^{-}_0(0)
+{\Delta}^{-}_0(0) +{\Delta }^{+}_0(0)$ & $2+4$ & $1+3$ \\
\end{tabular} 
\end{center}
\begin{center} 
Table 6. Coadjoint orbits of the Poincar\'{e} group ${\Oo ({\R }^6, K)}_{e_5}$. 
\end{center}

\section{Normal forms}\label{section 7}

We now give a table of explicit tuples $({\R }^6,Y,e_5;K)$ which
represent the corresponding affine cotypes listed in 
table 6. \medskip

In our list of normal forms we use the following conventions.
Let $\mathfrak{e} = \{e_0, e_1, \ldots , e_4, e_5 \} $ be the
{\sl standard basis} for ${\R }^6$ such that the Gram matrix of
the inner product is $K=${\tiny $\begin{pmatrix}
0 & 0 & 1 \\ 0 & G & 0 \\ 1 & 0 & 0 \end{pmatrix} $},
where $G =${\tiny $\begin{pmatrix} -I_3 & 0 \\ 0 & 1
\end{pmatrix} $}. We call $K$ the {\sl standard form}
of the inner product $\gamma $ on ${\R }^6$ and $G$ the
{\sl standard form} of the Lorentz inner product on ${\R }^4$ with
standard basis $\widetilde{\mathfrak{e}} = \{ e_1, \ldots , e_4 \} $.
\medskip

If $Y \in \oo ({\R }^6,K)$, then the matrix of $Y$ with respect to
the standard basis $\mathfrak{e}$ is
\begin{displaymath}
\mbox{{\footnotesize $\left( \begin{array}{ccr}
a & -x^TG & 0 \\
y & \widetilde{Y} & x \\
0 & -y^TG & -a \end{array} \right) $,}}
\end{displaymath}
where $a\in \R $, $x,y \in {\R }^4$ and $\widetilde{Y} \in 
\oo ({\R }^4, G)$.
Thus
with respect to the standard basis $\widetilde{\mathfrak{e}}$ the
matrix of $\widetilde{Y}$ is
\begin{displaymath}
\mbox{{\footnotesize $\left( \begin{array}{cc}
\widehat{z} & b \\ b^T & 0 \end{array} \right) $}} =
\mbox{{\footnotesize $\left( \begin{array}{rrrc}
0 & -z_3 & z_2 & b_1 \\
z_3 & 0 & -z_1 & b_2 \\
-z_2 & z_1 & 0 & b_3 \\
b_1 & b_2 & b_3 & 0 \end{array} \right) $},}
\end{displaymath}
where $b,z \in {\R }^3$. In other words,
\begin{displaymath}
Y= \mbox{{\footnotesize $\left( \begin{array}{c|rrrr|r}
a & x_1 & x_2 & x_3 & -x_4 & 0 \\ \hline
y_1 & 0 & -z_3 & z_2 & b_1 & x_1 \\
y_2 & z_3 & 0 & -z_1 & b_2 & x_2 \\
y_3 & -z_2 & z_1 & 0 & b_3 & x_3 \\
y_4 & b_1& b_2&b_3 & 0 & x_4 \\ \hline
0 & y_1 & y_2 & y_3 & -y_4 & -a \end{array} \right) $}.}
\end{displaymath}

In the list of normal forms below we give the matrix $\widetilde{Y}$ 
and vector $v$ of the little cotype that follows, we assume that 
the given the little cotype, represented by 
$({\R }^4, \widetilde{Y}, v; G)$. The normal form matrix $Y$ of 
the corresponding to the cotype represented by $({\R }^6, Y, e_5; K)$ is 
\begin{displaymath}
Y=\mbox{{\tiny $\left( \begin{array}{c|c|c} 
0 & -v^TG & 0 \\ \hline 
\verysmallrowspace 0 & \widetilde{Y} & v \\ \hline 
0 & 0 & 0 \end{array} \right) $}}. 
\end{displaymath}

Below is a list of representatives of the affine cotypes given in table 6.
\begin{itemize}
\item[1.] \textbf{Affine cotype}: ${\nabla }^{-}_5(0), \mu +
{\Delta }^{-}_0(0)$. 
\item[\mbox{}]{\sl sum basis}: $\{ {\mu }^{-2}Y^4w,\, -{\mu }^{-2}Y^3w,
\, {\mu }^{-1}Y^2w, \, Yw, \, w;  \, z \} $. {\sl conditions}: $Y^5w = 
Yz =0$; $\gamma (w,Y^4w) = -{\mu }^2$, $\gamma (z,z) =-1$. 
\item[\mbox{}] \textbf{little cotype}: Normal form basis:
\begin{displaymath}
\{ \ttfrac{1}{\sqrt{2}} ({\mu }^{-2}Y^3w - Yw),
\, {\mu }^{-1}Y^2w, \, z, \, \ttfrac{1}{\sqrt{2}}
({\mu }^{-2}Y^3w +Yw) \} . 
\end{displaymath}
Normal form matrix $\widetilde{Y}$ and vector $v$. 
\begin{displaymath}
\widetilde{Y} = 
\mbox{{\tiny $\left( \begin{array}{rc}
-\ttfrac{\mu }{\sqrt{2}} \, \widehat{e_3} &
\ttfrac{\mu }{\sqrt{2}} \, e_2 \\
\smallrowspace \ttfrac{\mu }{\sqrt{2}} \, e^T_2 & 0 \end{array} \right) $};} 
\, \, \, v = \ttfrac{1}{\sqrt{2}}(-e_1+e_4). 
\end{displaymath}
\end{itemize}
\begin{itemize}
\item[2.] \textbf{Affine cotype}: ${\nabla }_4(0,0) +
{\Delta }^{-}_0(i\beta , \mathrm{IP})$. 
\item[\mbox{}]{\sl sum basis}: $\{ -z ,\, Yz,\, Yw, w;\, u,
\, {\beta }^{-1} Yu \} $. {\sl conditions}: $Y^2w = Y^2z =0$,
$(Y^2 +{\beta }^2)u=0$; $\gamma (Yz, w) = 1$ and $\gamma (u,u) =-1$. 
\item[\mbox{}] \textbf{little cotype}: Normal form basis: 
\begin{displaymath}
\{  
\ttfrac{1}{\sqrt{2}} (Yw+z),
\, u, \, {\beta }^{-1}Yu, \, \ttfrac{1}{\sqrt{2}}
(Yw-z)  \} .  
\end{displaymath} 
Normal form matrix $\widetilde{Y}$ 
and vector $v$: 
\begin{displaymath}
\widetilde{Y} = 
\mbox{{\tiny $\left( \begin{array}{cc}
\beta \, \widehat{e_1} & 0 \\
\smallrowspace 0 & 0 \end{array} \right) $};} \, \, \, 
v = \ttfrac{1}{\sqrt{2}}(e_1+e_4) .
\end{displaymath}
\end{itemize}
\begin{itemize}
\item[3.] \textbf{Affine cotype}: ${\nabla }_4(0,0) +
{\Delta }^{-}_0(0) + {\Delta }^{-}_0(0)$. 
\item[\mbox{}]{\sl sum basis}: $\{ -z,\, Yz, \, Yw,\,  w;\, u; \, v \} $. 
{\sl conditions}: $Y^2w = Y^2z = Yu =Yv =0$; 
$\gamma (Yz, w) = 1$, $\gamma (u,u) = \gamma (v,v) =-1$. 
\item[\mbox{}] \textbf{little cotype}: Normal form basis: 
$ \{ \ttfrac{1}{\sqrt{2}} (z + Yw),u, \, v, \, \ttfrac{1}{\sqrt{2}}(Yw-z) \} $.
Normal form matrix $\widetilde{Y}$ and vector $v$: $\widetilde{Y} =  
0; \, \, \, v = \ttfrac{1}{\sqrt{2}}(e_1+e_4) $.
\end{itemize}
\begin{itemize}
\item[4.] \textbf{Affine cotype}: ${\nabla }^{-}_3(0), \mu +
{\Delta }^{+}_2(0) $. 
\item[\mbox{}]{\sl sum basis}: $\{ {\mu }^{-2}Y^2w, \, {\mu }^{-1}Yw ,\, w;
\, u, \, Yu, \, Y^2u \} $. {\sl conditions}: $Y^3w = Y^3u =0$; 
$\gamma (w, Y^2w) = {\mu }^2$, $\gamma (u,Y^2u) = 1$. 
\item[\mbox{}] \textbf{little cotype}: Normal form basis:
\begin{displaymath}
\{ {\mu }^{-1}Yw, \,
\ttfrac{1}{\sqrt{2}} (u - Y^2u),
\, Yu, \, \ttfrac{1}{\sqrt{2}}(u +Y^2u)  \} .
\end{displaymath} 
Normal form matrix $\widetilde{Y}$ 
and vector $v$: $\widetilde{Y} = 
\mbox{{\tiny $\begin{pmatrix}
\ttfrac{1}{\sqrt{2}}\, \widehat{e_1} & \ttfrac{1}{\sqrt{2}}\, e_3 \\
\smallrowspace \ttfrac{1}{\sqrt{2}}\, e^T_3 & 0 \end{pmatrix} $};} 
\, \, \, v = \mu \, e_1 $.
\end{itemize}
\begin{itemize}
\item[5] \textbf{Affine cotype}: ${\nabla }^{-}_3(0), \mu  +
{\Delta }^{-}_0(i\beta , \mathrm{IP}) +{\Delta }^{+}_0(0) $. 
\item[\mbox{}]{\sl sum basis}: $\{ {\mu }^{-2}Y^2w, \,
{\mu }^{-1}Yw ,\, w;
\, u, \, {\beta }^{-1}Yu; \, v \} $. {\sl conditions}: $Y^3w = Yv =0$,
and $(Y^2+{\beta }^2)u=0$; $\gamma (w, Y^2w) = {\mu }^2$, $\gamma (u,u) = -1$,
and $\gamma (v,v) =1$. 
\item[\mbox{}] \textbf{little cotype}: Normal form basis: $\{  {\mu }^{-1}Yw,
\, u, \, {\beta }^{-1}Yu, \, v \} $. Normal form matrix $\widetilde{Y}$ and 
vector $v$: $\widetilde{Y} = 
\mbox{{\tiny $\begin{pmatrix}
\beta\, \widehat{e_1} & 0 \\
\smallrowspace 0 & 0 \end{pmatrix} $};} \, \, \, v = \mu \, e_1$.
\end{itemize}
\begin{itemize}
\item[6.] \textbf{Affine cotype}: ${\nabla }^{-}_3(0), \mu  +
{\Delta }_0(\alpha , \mathrm{RP}) +{\Delta }^{-}_0(0) $. 
\item[\mbox{}]{\sl sum basis}: $\{ {\mu }^{-2}Y^2w, \, {\mu }^{-1}Yw ,\, w;
\, u, \, {\alpha }^{-1}Yu; \, v \} $. {\sl conditions}: $Y^3w = Yv =0$,
and $(Y^2-{\alpha }^2)u=0$; $\gamma (w, Y^2w) = {\mu }^2$, $\gamma (u,u) = 1$,
and $\gamma (v,v) =-1$. 
\item[\mbox{}] \textbf{little cotype}: Normal form basis: $\{ {\mu }^{-1}Yw, 
\, {\alpha }^{-1}Yu, \, v, \, u;\, w \} $. Normal form matrix 
$\widetilde{Y}$ and vector $v$: $\widetilde{Y} = 
\mbox{{\tiny $\begin{pmatrix}
0 & \alpha \, e_2  \\
\rule{0pt}{7pt} \alpha \, e^T_2 & 0 \end{pmatrix} $};} \, \, \, 
v = \mu \, e_1 $.
\end{itemize}
\begin{itemize}
\item[7.] \textbf{Affine cotype}: ${\nabla }^{-}_3(0), \mu  +
{\Delta }^{-}_0(0) +{\Delta }^{-}_0(0) +{\Delta }^{+}_0(0) $. 
\item[\mbox{}]{\sl sum basis}: $\{ {\mu }^{-2}Y^2w, \, {\mu }^{-1}Yw ,\, w;
\, u; \, v; \, z \} $. {\sl conditions}: $Y^3w = Yu = Yv = Yz =0$; 
$\gamma (w, Y^2w) = {\mu }^2$, $\gamma (u,u) =
\gamma (v,v) = -1$, $\gamma (z,z) =1$.
\item[\mbox{}] \textbf{little cotype}: Normal form basis: $\{ {\mu }^{-1}Yw, 
\, u, \, v, \, z \} $. Normal form matrix $\widetilde{Y}$ and 
vector $v$: $\widetilde{Y} =  0; \, \, \, v = \mu \, e_1$.
\end{itemize}
\begin{itemize}
\item[8.] \textbf{Affine cotype}: ${\nabla }^{+}_3(0), \mu  +
{\Delta }^{-}_0(i\beta , \mathrm{IP}) +{\Delta }^{-}_0(0) $. 
\item[\mbox{}]{\sl sum basis}: $\{ -{\mu }^{-2}Y^2w, 
\, {\mu }^{-1}Yw ,\, w;\,
u, \, {\beta }^{-1}Yu; \, v \} $. {\sl conditions}: $Y^3w = Yv =0$ and
$(Y^2+{\beta }^2)u =0$; $\gamma (w, Y^2w) = -{\mu }^2$, and $\gamma (u,u) =
\gamma (v,v) =-1$.
\item[\mbox{}] \textbf{little cotype}: Normal form basis: $\{ u, 
{\beta }^{-1}Yu, \, v, \, {\mu }^{-1}Yw \} $. Normal form 
matrix $\widetilde{Y}$ and vector $v$: $\widetilde{Y} = 
\mbox{{\tiny $\begin{pmatrix}
\beta\, \widehat{e_3} & 0 \\
\smallrowspace 0 & 0 \end{pmatrix} $};} \, \, \, v = \mu \, e_4$.
\end{itemize}
\begin{itemize}
\item[9.] \textbf{Affine cotype}: ${\nabla }^{+}_3(0), \mu  +
{\Delta }^{-}_0(0) +{\Delta}^{-}_0(0) + {\Delta }^{-}_0(0)$. 
\item[\mbox{}]{\sl sum basis}: $\{ -{\mu }^{-2}Y^2w, \, {\mu }^{-1}Yw ,\, w;
\, u;  \, v; \, z \} $. {\sl conditions}: $Y^3w = Yu = Yv = Yz =0$; 
$\gamma (w, Y^2w) = -{\mu }^2$, and $\gamma (u,u) =
\gamma (v,v) = \gamma (z,z) =-1$.
\item[\mbox{}] \textbf{little cotype}: Normal form basis: $\{ u,\, v,
\, z, \, {\mu }^{-1}Yw \} $. Normal form matrix $\widetilde{Y}$ 
and vector $v$: $\widetilde{Y} =  0; \, \, \, v = \mu \, e_4$.
\end{itemize}
\begin{itemize}
\item[10.] \textbf{Affine cotype}: ${\nabla }_2(0,0) +
{\Delta }^{+}_2(0) +{\Delta}^{-}_0(0)$. 
\item[\mbox{}]{\sl sum basis}: $ \{ z,\, w;\, Y^2u, \, Yu ,\, u;\, v \} $. 
{\sl conditions}: $Y^3u = Yw = Yv = Yz =0$; $\gamma (z,w) = 1$, 
$\gamma (u,Y^2u) = 1$, and
$\gamma (v,v) =-1$.
\item[\mbox{}] \textbf{little cotype}: Normal form basis: 
\begin{displaymath}
\{ \ttfrac{1}{\sqrt{2}}(u-Y^2u),
\, Yu, \, v,\, \ttfrac{1}{\sqrt{2}}(u+Y^2u) \} .
\end{displaymath}
Normal form matrix $\widetilde{Y}$ and 
vector $v$: $\widetilde{Y} = 
\mbox{{\tiny $\begin{pmatrix}
\ttfrac{1}{\sqrt{2}}\, \widehat{e_3} & \ttfrac{1}{\sqrt{2}}\, e_2  \\
\smallrowspace \ttfrac{1}{\sqrt{2}}\, e^T_2 & 0 \end{pmatrix} $};} \, \, 
\, v = 0$.
\end{itemize}
\begin{itemize}
\item[11.] \textbf{Affine cotype}: ${\nabla }_2(0,0) +
{\Delta }^{-}_0(i\beta , \mathrm{IP}) +{\Delta}_0(\alpha , \mathrm{RP})$. 
\item[\mbox{}]{\sl sum basis}: $\{ z,\, w;\, u,\, {\beta }^{-1} Yu ;\, v, \,
{\alpha }^{-1}Yv  \} $. {\sl conditions}: $Yz = Yw =0$, $(Y^2+{\beta }^2)u =0$,
$(Y^2-{\alpha }^2)v =0$; $\gamma (z,w) = \gamma (v,v) = 1$, and
$\gamma (u,u) =-1$.
\item[\mbox{}] \textbf{little cotype}: Normal form basis: $\{  u, 
\, {\beta }^{-1}Yu, \, {\alpha }^{-1}Yv, \, v \} $.
Normal form matrix $\widetilde{Y}$ and 
vector $v$: $\widetilde{Y} = 
\mbox{{\tiny $\begin{pmatrix}
\beta \, \widehat{e_3} & \alpha \, e_3  \\
\rule{0pt}{7pt} \alpha \, e^T_3 & 0 \end{pmatrix} $};} \, \, \, v=0 $.
\end{itemize}
\begin{itemize}
\item[12.] \textbf{Affine cotype}: ${\nabla }_2(0,0) +
{\Delta }^{-}_0(i\beta , \mathrm{IP}) +{\Delta}^{-}_0(0)+{\Delta}^{+}_0(0)$. 
\item[\mbox{}]{\sl sum basis}: $\{ z,\, w;\, u,\, {\beta }^{-1} Yu ;\, v;
\, y  \} $. {\sl conditions}: $Yz = Yw =Yv= Yy =0$; 
$\gamma (z,w) = \gamma (y,y) = 1$ and
$\gamma (u,u) = \gamma (v,v) =-1$.
\item[\mbox{}] \textbf{little cotype}: Normal form basis: $\{ u,
\, {\beta }^{-1}Yu, \, v;\, y; \, w \} $. Normal form matrix 
$\widetilde{Y}$ and 
vector $v$: $\widetilde{Y} = 
\mbox{{\tiny $\begin{pmatrix}
\beta \, \widehat{e_3} & 0  \\
\rule{0pt}{7pt} 0 & 0 \end{pmatrix} $};} \, \, \, v=0 $.
\end{itemize}
\begin{itemize}
\item[13.] \textbf{Affine cotype}: ${\nabla }_2(0,0) +
{\Delta }_0(\alpha , \mathrm{RP}) +{\Delta}^{-}_0(0)+{\Delta}^{-}_0(0)$. 
\item[\mbox{}]{\sl sum basis}: $\{ z,\, w;\, u,\, {\alpha }^{-1} Yu ;
\, v; \, y \} $. {\sl conditions}: $Yz = Yw =Yv= Yy =0$,
$(Y^2 - {\alpha }^2)u =0$; $\gamma (z,w) = \gamma (u,u) =1$ and
$ \gamma (v,v) = \gamma (y,y) =-1$.
\item[\mbox{}] \textbf{little cotype}: Normal form basis: $\{ 
{\alpha }^{-1}Yu, \, v, \, y, \, u  \} $. 
Normal form matrix $\widetilde{Y}$ and vector $v$: $\widetilde{Y} = 
\mbox{{\tiny  $\begin{pmatrix}
0 & \alpha \, e_1  \\
\rule{0pt}{7pt} \alpha \, e^T_1 & 0 \end{pmatrix} $};} \, \, \, v=0$.
\end{itemize}
\begin{itemize}
\item[14.] \textbf{Affine cotype}: ${\nabla }_2(0,0) +
{\Delta }^{-}_0(0) +{\Delta}^{-}_0(0)+{\Delta}^{+}_0(0)$. 
\item[\mbox{}]{\sl sum basis}: $ \{ z,\, w;\, u ;\, v; \, y;\, s  \} $. 
{\sl conditions}: $Yz = Yw =Yv= Yy =Ys = 0$; $\gamma (z,w) = \gamma(s,s)= 1$ 
and
$\gamma (u,u) = \gamma (v,v) = \gamma (y,y) =-1$.
\item[\mbox{}] \textbf{little cotype}: Normal form basis: $\{ v, 
\, u, \, y, \, s \} $. Normal form matrix $\widetilde{Y}$ and vector $v$: 
$\widetilde{Y} =  0; \, \, \, v =0$.
\end{itemize}

\section{Acknowledgements}
We thank Hans Duistermaat for useful discussions.

\end{document}